\documentclass{amsart}
\usepackage{latexsym}
\usepackage{amssymb}
\usepackage{amsmath}
\usepackage{tikz}
\usepackage[all]{xy}
\usepackage{url}
\def\ni{\noindent}
\def\be{\begin{equation}}
\def\ee{\end{equation}}
\def\bC{\mathbb{C}}
\def\bR{\mathbb{R}}
\def\bH{\mathbb{H}}
\def\gg{\mathfrak{g}}
\def\gk{\mathfrak{k}}
\def\gp{\mathfrak{p}}
\def\gt{\mathfrak{t}}
\def\D{\mathcal{D}}
\newtheorem{thm}{Theorem}[section]

\newtheorem{lem}[thm]{Lemma}
\newtheorem{prop}[thm]{Proposition}
\theoremstyle{definition}

\theoremstyle{remark}
\newtheorem{rem}[thm]{Remark}
\numberwithin{equation}{section}
\begin{document}

\title[]{On a characteristic of the first eigenvalue of the Dirac operator
 on compact spin symmetric spaces with a K\"ahler or Quaternion-K\"ahler structure}%
\author{Jean-Louis Milhorat}%
\address{Laboratoire Jean Leray, UMR 6629,
Universit\'e de Nantes}%
\email{jean-louis.milhorat@univ-nantes.fr}%

%\thanks{}%
%\keywords{}%

%\date{}%
%\dedicatory{}%
%\commby{}%
% ----------------------------------------------------------------
\begin{abstract}
It is shown that on a compact spin symmetric space with a K\"ahler or Quaternion-K\"ahler structure, the
first eigenvalue of the Dirac operator is linked to a ``{lowest}'' action of the holonomy, given by the fiberwise action
on spinors of the canonical forms characterized by this holonomy. The result is also verified for the symmetric space
$\mathrm{F}_4/\mathrm{Spin}_9$, proving that it is valid for all the ``{possible}'' holonomies in the Berger's list
occurring in that context. The proof is based on a characterization of the first eigenvalue of the Dirac operator given in \cite{Mil05} and \cite{Mil06}.
By the way, we review an incorrect statement in the proof of the first lemma in \cite{Mil05}.
\end{abstract}
\maketitle
% ----------------------------------------------------------------
\section{Introduction}
Let $(M^n,g)$ be a spin compact Riemannian manifold with positive scalar curvature, more precisely such that
$\mathrm{Scal}_{\circ}:=\min_{m\in M}\mathrm{Scal}(m)
>0$. Under this assumption, the only groups $G$ in Berger's list such that the restricted holonomy group of $M$
verifies $Hol^{\circ}\subset G$ are (cf. for instance \cite{Bes}) $G=U_m$, $n=2m$, ($M$ is then K\"ahler), $G=\mathrm{Sp}_{m}\cdot \mathrm{Sp}_1$,
 $n=4m$, ($M$ is then Quaternion-K\"ahler) or $G=\mathrm{Spin}_9$, $n=16$, and then $M$ is isometric to the Cayley
 plane $\mathbb{O}P^2=\mathrm{F}_4/\mathrm{Spin}_9$, \cite{Ale}, \cite{BG72}.

\noindent Assuming $n=4m$ in order to compare all the possible cases, there exist sharp lower bounds for the
square of the eigenvalues of the Dirac operator whose dependence on the holonomy is summarized
in the following illustration:
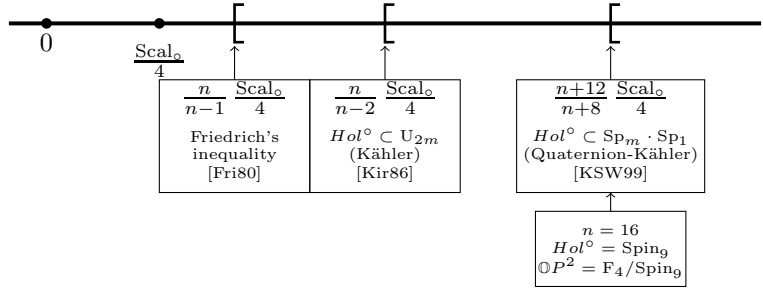
\begin{center}
\begin{tikzpicture}
 \draw [line width=1.5pt] (-5,0) -- (5,0);\fill (-4.5,0) circle(2pt);
 \draw (-4.5,-0.25) node {$0$};\fill (-3,0) circle(2pt);
 \draw (-3,-0.5) node {$\frac{\mathrm{Scal}_{\circ}}{4}$};
 \draw [line width=1pt](-1.875,0.25)--(-2,0.25)--(-2,-0.25)--(-1.875,-0.25);
 \draw [black] (-3,-2.25) rectangle (-1,-0.75);
 \draw [->] (-2,-0.75)--(-2,-0.3);
 \draw (-2,-1) node {$\frac{n}{n-1}\,\frac{\mathrm{Scal}_{\circ}}{4}$};
 \draw (-2, -1.5) node {\tiny{Friedrich's}};
 \draw (-2,-1.75) node {\tiny{inequality}};
 \draw (-2,-2) node {\tiny{\cite{Fri80}}};
 \draw [line width=1pt](0.125,0.25)--(0,0.25)--(0,-0.25)--(0.125,-0.25);
 \draw [black] (-1,-2.25) rectangle (1,-0.75);
 \draw [->] (0,-0.75)--(0,-0.3);
 \draw (0,-1) node {$\frac{n}{n-2}\,\frac{\mathrm{Scal}_{\circ}}{4}$};
 \draw (0, -1.5) node {\tiny{$Hol^{\circ}\subset \mathrm{U}_{2m}$}};
 \draw (0,-1.75) node {\tiny{(K\"ahler)}};
 \draw (0,-2) node {\tiny{\cite{Kir86}}};
 \draw [line width=1pt](3.125,0.25)--(3,0.25)--(3,-0.25)--(3.125,-0.25);
 \draw [black] (1.75,-2.25) rectangle (4.25,-0.75);
 \draw [->] (3,-0.75)--(3,-0.3);
 \draw (3,-1) node {$\frac{n+12}{n+8}\,\frac{\mathrm{Scal}_{\circ}}{4}$};
 \draw (3, -1.5) node {\tiny{$Hol^{\circ}\subset \mathrm{Sp}_{m}\cdot \mathrm{Sp}_1$}};
 \draw (3, -1.75) node {\tiny{(Quaternion-K\"ahler)}};
 \draw (3,-2) node {\tiny{\cite{KSW99}}};
\draw (3,-2.75) node {\tiny{$n=16$}};
 \draw (3,-3) node {\tiny{$Hol^{\circ}= \mathrm{Spin}_9$}};
 \draw (3,-3.25) node {\tiny{$\mathbb{O}P^2=\mathrm{F}_4/\mathrm{Spin}_9$}};
 \draw [black] (2,-3.5) rectangle (4,-2.5);
 \draw [->] (3,-2.5)--(3,-2.25);
\end{tikzpicture}
\end{center}
The study of limiting manifolds, that are manifolds for which there exists a spinorfield $\Psi$ such that
\be\label{lim.psi}\D^2\Psi=\lambda^2\, \Psi\,,\ee where $\D$ is the Dirac operator, and where $\lambda^2$
is one of the bounds quoted above is due to C. B\"ar in the general
case, \cite{Bar93}, to A. Moroianu in the case of K\"ahler manifolds, \cite{Mor95}, \cite{Mor99},
and to W. Kramer, U. Semmelmann and G. Weingart in the case of Quaternion-K\"ahler manifolds, \cite{KSW98a}.

\ni The study of limiting manifolds in the K\"ahler and Quaternion-K\"ahler cases involves a special
condition for spinor fields $\Psi$ verifying \eqref{lim.psi}, which is linked to the decomposition of
the spinor space $\Sigma$ into irreducible components under the action of the holonomy group.

\ni The spinor bundle $\Sigma(M)$ of a spin K\"ahler manifold $(M,g,J)$ of complex dimension $m$ decomposes into
a sum of parallel subbundles $\Sigma(M)=\oplus_{r=0}^{m}\Sigma_{r}(M)$, each section of the bundle
$\Sigma_{r}(M)$ being an eigenvector for the eigenvalue $i\,(m-2r)$, for the fiberwise action of the K\"ahler form
$\Omega$ on spinors, \cite{Kir86}.

\ni It is then a characteristic of limiting K\"ahler manifolds that in the space of spinor fields $\Psi$
verifying \eqref{lim.psi}, there always exists an element such that
$\Omega\cdot \Psi =0$ if $m$ is even, or $\Omega\cdot \Psi =\pm i$, if $m$ is odd. We also may formulate this property as~:

$\bullet$ In the space of spinor fields $\Psi$ verifying \eqref{lim.psi}, there always exists an element such that
\be\label{low.hol.1}\frac{\|\Omega\cdot \Psi\|^2}{\|\Psi\|^2}\quad \text{ is minimal.}
 \ee

\ni In order to illustrate how this property corresponds to a ``{lower action}'' of the ``{K\"ahler holonomy}'',
recall that the above decomposition of the spinor bundle corresponds to the decomposition of the spinor space
$\Sigma_{2m}$ under the action of the groups $\mathrm{U}_1\times \mathrm{SU}_m$ when $m$ is even, or
$\mathrm{S}({U}_1\times {U}_m)$ when $m$ is odd, actions given by the commutative diagrams

$$\xymatrix{ \mathrm{U}_1\times \mathrm{SU}_m \ar[r]\ar[rd]^{(m\,\mathrm{even})} & \mathrm{Spin}_{2m}\ar[d]&&
\mathrm{S}(\mathrm{U}_1\times \mathrm{U}_m) \ar[r]\ar[rd]^{(m\,\mathrm{odd})} & \mathrm{Spin}_{2m}\ar[d]\\
& \mathrm{U}_{m}\subset SO_{2m}&&& \mathrm{U}_{m}\subset SO_{2m}
}$$

This decomposition may be easily expressed in the case $m$ even, \cite{Sal92},
\be\label{decomp.Kahl} \Sigma_{2m} = \bigoplus_{r=0}^{m} L^{-\frac{1}{2}(m-2r)}\otimes \Lambda^{r} E\,,
\ee
where $L^{k}:=L^{\otimes k}$ is the $k$-symmetric power of the standard representation $L$ of $\mathrm{U}_1$,
and $E$ is the standard representation of the group $\mathrm{SU}_m$.

\ni The action of the K\"ahler form corresponds to twice the action of the Lie algebra of $\mathrm{U}_1$ on $\Sigma_{2m}$,
hence is equal to $i\, (m-2r)\, \mathrm{id}$ on each component of the sum
\eqref{decomp.Kahl}. So, limiting K\"ahler manifolds of even complex dimension are characterized by the existence
of a spinor field verifying \eqref{lim.psi}, which is also a section of the bundle corresponding to the component
$L^{0}\otimes \Lambda^{m/2}$ in \eqref{decomp.Kahl}, this component being, roughly speaking, the one with
minimal ``{$\mathrm{U}_1$ holonomy}'' in \eqref{decomp.Kahl}. There is an analogous description in the case $m$ odd.

\ni There exists an analogous criterion for Quaternion-K\"ahler limiting manifolds.
The spinor bundle $\Sigma(M)$ of a spin Quaternion-K\"ahler manifold $(M,g)$ of dimension $4m$ decomposes into
a sum of parallel subbundles $\Sigma(M)=\oplus_{r=0}^{m}\Sigma_{r}(M)$ , each section of the bundle
$\Sigma_{r}(M)$ being an eigenvector for the eigenvalue $6m-4r(r+2)$, for the fiberwise action of the fundamental
$4$-form $\Omega$ on spinors, \cite{HM95a}.
This decomposition corresponds to the decomposition of the spinor space $\Sigma_{4m}$
into irreducible components under the action of the group $\mathrm{Sp}_m\times \mathrm{Sp}_1$
given by the commutative diagram
$$\xymatrix{ \mathrm{Sp}_m\times \mathrm{Sp}_1 \ar[r]\ar[rd] & \mathrm{Spin}_{4m}\ar[d]\\
& \mathrm{Sp}_{m}.\mathrm{Sp}_1\subset SO_{4m}
}$$
One gets \cite{BS83}, \cite{Wan89}, \cite{HM95a},
\be\label{decomp.QK} \Sigma_{4m} = \bigoplus_{r=0}^{m} \mathrm{S}^{r}H\otimes \Lambda_{\circ}^{m-r} E\,,
\ee
where $\mathrm{S}^{k} H$ is the $k$-symmetric power of the standard representation $H$ of $\mathrm{Sp}_1$ in
the space $\bH\simeq\bC^2$, $E$ is the standard representation of the group $\mathrm{Sp}_m$ in the space
$\bH^{m}\simeq\bC^{2m}$, and $\Lambda_{\circ}^{k} E$ is the irreducible hermitian complement of
$\omega\wedge \Lambda_{\circ}^{k-2} E$
in $\Lambda^{k} E$, $\omega$ being the standard symplectic form on $E$.

\ni Quaternion-K\"ahler limiting manifolds are characterized by the fact that,
 among the spinor fields $\Psi$ verifying \eqref{lim.psi},
there always exists a section of the bundle corresponding to the space
$\mathrm{S}^{0}H\otimes \Lambda_{\circ}^{m} E$.

\ni This space may be characterized by the fact that it is the only component in the decomposition
\eqref{decomp.QK} on which the Casimir operator of the subgroup $\mathrm{Sp}_1$ acts trivially (roughly
speaking, one could say that this is the component with minimal ``{$\mathrm{Sp_1}$ holonomy}'' in \eqref{decomp.QK}).
From a geometrical point of view,
the invariant given by this action of the Casimir operator, gives raise to the operator $\Omega -6m\,
\mathrm{id}$, where $\Omega$ is the $4$-fundamental form $\Omega$ acting on spinor fields by Clifford multiplication.
Hence, we may re-formulate the characterization of Quaternion-K\"ahler limiting manifolds as:
\begin{itemize}
\item In the space of spinor fields $\Psi$ verifying \eqref{lim.psi}, there always exists an element such that
$(\Omega-6m\, \mathrm{id})\cdot \Psi =0$.
\item In the space of spinor fields $\Psi$ verifying \eqref{lim.psi}, there always exists an element such that
\be\label{low.hol.2}\frac{\|(\Omega-6m\,\mathrm{id})\cdot \Psi\|^2}{\|\Psi\|^2}\quad \text{ is minimal.}
 \ee
\end{itemize}

\ni Now, the following example makes think that the above criteria are not characteristic of limiting manifolds.
Consider the grassmannian  $\mathrm{Gr}_2(\mathbb{C}^{m+2})=\frac{\mathrm{SU}_{m+2}}{\mathrm{S}(\mathrm{U}_2\times
\mathrm{U}_m)}$, endowed with its canonical metric induced by the Killing form of $\mathrm{SU}_{m+2}$
sign-changed, which is both K\"ahler and Quaternion-K\"ahler (cf. 14.53 in \cite{Bes} for details).
It is shown in \cite{Mil98} that there exists a spinorfield $\Psi$ such that \eqref{lim.psi} is verified for the
first eigenvalue $\lambda$ and $\eqref{low.hol.2}$ is verified for the fundamental ``{Quaternion-K\"ahler}''
$4$-form. Actually, it may also be checked that $\eqref{low.hol.1}$ is also verified for the K\"ahler form.

\ni Hence it seems natural to conjecture that the above property is not a characteristic of limiting manifolds and
the aim of this paper is to prove that the conjecture is true for compact spin symmetric spaces.

\begin{prop}\label{result} Any
 spin compact simply connected irreducible symmetric space $G/K$ of ``type I'',
endowed with a K\"ahler or Quaternion-K\"ahler structure, verifies the following criterion.

\ni  Among the spinorfields $\Psi$ verifying $\D^2\Psi=\lambda^2\, \Psi$,
where $\lambda$ is the first eigenvalue of the Dirac operator, there exists at least one such that
$$\frac{\|\Omega\cdot \Psi\|^2}{\|\Psi\|^2}\,,\quad\text{respectively}\quad
\frac{\|(\Omega-6m\,\mathrm{id})\cdot \Psi\|^2}{\|\Psi\|^2}\,,\quad\text{is minimal,}$$
where $\Omega$ is respectively the K\"ahler form or the fundamental ``{Quaternion-K\"ahler}'' $4$-form
of the manifold under consideration.

\ni There is an analogous result for the Cayley plane $\mathbb{O}P^2=\mathrm{F}_4/\mathrm{Spin}_9$,
$\Omega$ being in this case the canonical $8$-form on
manifolds with holonomy $\mathrm{Spin}_9$.
\end{prop}

\section{Preliminaries for the proof}

\subsection{Spectrum of the Dirac operator on spin compact irreducible symmetric spaces}
We consider a spin compact simply connected irreducible symmetric space $G/K$ of ``type I'',
where $G$ is a simple compact and simply-connected Lie group and $K$ is the
connected subgroup formed by the fixed elements of an involution
$\sigma$ of $G$.  This involution induces the
Cartan decomposition of the Lie algebra $\gg$ of $G$ into
$$\gg=\gk\oplus\gp\,,$$ where $\gk$ is the Lie algebra of $K$ and
$\gp$ is the vector space $\{X\in\gg\,;\, \sigma_{*}\cdot X=-X\}$.
The symmetric space $G/K$ is endowed with the Riemannian metric induced by the restriction to $\gp$ of
the Killing form of $G$ sign-changed.

\ni The spin condition implies that the homomorphism
$$\alpha: h\in K \longmapsto \mathrm{Ad}_{G}(h)_{|\gp}\in \mathrm{SO}(\gp)$$
lifts to a homomorphism $\widetilde{\alpha}:H\longrightarrow\mathrm{Spin}(\gp)$
such that $\xi\circ \widetilde{\alpha}=\alpha$ where $\xi$ is the two-fold covering
$\mathrm{Spin}(\gp)\rightarrow \mathrm{SO}(\gp)$, \cite{CG}.

\ni Then the group $K$ inherits a spin representation given by
$$\widetilde{\rho_K}: K \xrightarrow{\widetilde{\alpha}} \mathrm{Spin}(\gp)\xrightarrow{\rho}
\mathrm{GL}_{\bC}{(\Sigma)}\,,$$
where $\rho$ is the spinor representation in the complex spinor space $\Sigma$.

\ni The Dirac operator has a real discrete spectrum, symmetric with respect to the origin.
A real number $\lambda$ belongs to the spectrum if and only if there exists an irreducible representation
$\gamma:G\rightarrow \mathrm{GL}_{\bC}{(V_{\gamma})}$ whose restriction $\mathrm{Res}_K^G(\gamma)$ to the subgroup $K$, contains in its decomposition
into irreducible parts, a representation equivalent to some irreducible component of the decomposition
 of the spin representation $\widetilde{\rho_K}$ of $K$.
Then
\be\label{eigenv}\lambda^2=c_{\gamma}+n/16\,,\ee
where $c_\gamma$ is the Casimir eigenvalue of the irreducible representation $\gamma$ (which only depens on the equivalence class of $\gamma$)
and where $n=\dim (G/K)$, $n/16$ being $\mathrm{Scal}/8$ for the choice of the metric (cf. \cite{BH3M} or
 \cite{Gin09} for details).

\subsection{A characterization of the decomposition of the spin representation of $K$}

\ni We henceforward assume that
$G$ and $K$ have same rank and consider a fixed common maximal torus $T$. Hence $G/K$ has
even dimension $n=2m$.

\ni It is shown in \cite{Mil05} that the lowest eigenvalue $\lambda_{min}$ of the Dirac operator
verifies
\be\label{low.eig}\lambda_{min}^{2}=2\, \mathop{\min}_{1\leq i\leq
N}\|\beta_{i}\|^{2}+\frac{n}{8}\,,\ee
where the $\beta_{i}$, $1\leq i\leq N$, are the dominant weights (relative to the choice of $T$)
 occurring in the decomposition of the spin representation of $K$, and where the norm $\|\cdot\|$  on the space of weights
is induced by the Killing form of $G$ sign-changed.

\ni The proof of proposition~\eqref{result} is based on a characterization of those dominant weigths $\beta_{i}$ for which
$\|\beta_{i}\|^{2}$ is minimum, a characterization which appears implicitly in \cite{Mil06}. It is based on a
lemma of R. Parthasaraty in \cite{Par} (cf.~lemma~2.2), which gives the following characterization of dominant
weights occurring in the decomposition of the spin representation of $K$.

\ni Let $\Phi$ be the set of non-zero roots of the group $G$
 with respect to $T$. According to a classical terminology, a root
 $\theta$ is called compact if the corresponding root space is contained
 in $\gk_{\mathbb{C}}$ (that is, $\theta$ is a root of $K$ with respect to $T$) and
 non compact if the root space is contained in $\gp_{\mathbb{C}}$.
Let $\Phi_{G}^{+}$  be the set of positive roots of $G$,
$\Phi_{K}^{+}$
 be the set of positive roots of $K$,
and $\Phi_{\gp}^{+}$ be the set of positive non compact roots with
respect to a fixed lexicographic ordering in $\Phi$. The half-sums
of the positive roots of $G$ and $K$ are respectively denoted
$\delta_{G}$ and $\delta_{K}$ and the half-sum of non  compact
positive roots is denoted by $\delta_{\gp}$. The Weyl group of $G$
is denoted $W_{G}$. The space of weights is endowed with the
$W_{G}$-invariant scalar product $\langle,\rangle$ induced by the Killing
form of $G$ sign-changed.

\noindent To introduce the result of Parthasaraty first note that
 the common torus $T$ may be chosen in such a way that the weights of the spin representation of $K$
are
\begin{equation}\label{weights}
 \frac{1}{2}\, (\pm\alpha_1\pm\alpha_2\,\cdots\pm\alpha_m)\,,
\end{equation}
where $\alpha_1,\ldots,\alpha_m$ is an enumeration of the non compact positive roots,
the weights of the half-spin representations $\Sigma^{\pm}$ corresponding to
an even (resp. odd) number of negative signs.

\ni Thus weights of the spin representation of $K$ have the form

\be\label{weights.spin}
 \delta_{\gp} -\sum_{i\in I}\alpha_{i}\,,\qquad I\subset\{1,\ldots, m\}\,.
\ee

\begin{lem}[R. Parathasaraty, \cite{Par}] \label{Parth} Let
\begin{equation}W:=\{w\in W_{G}\;;\;
w\cdot \Phi_{G}^{+}\supset \Phi_{K}^{+}\}\,.
\end{equation}

\ni The spin representation of $K$ decomposes into irreducible components
as
\begin{equation}
\widetilde{\rho_{K}}=\bigoplus_{w\in W}\widetilde{\rho_{K}}_{|w}\;,
\end{equation}
where $\widetilde{\rho_{K}}_{|w}$ has for dominant weight
\begin{equation}
\beta_{w}:= w\cdot \delta_{G}-\delta_{K}\,.
\end{equation}
\end{lem}

\subsection{A characterization of highest weights of the spin representation of $K$ with minimal norm}

\begin{prop} \label{high.w.char} Let
\begin{align*}I_{-}:=&\{i\in\{1,\ldots,m\}\,,\, \langle \delta_K,\alpha_i\rangle< 0\}\,,
\intertext{and}
 I_0 :=&\{i\in\{1,\ldots,m\}\,,\, \langle \delta_K,\alpha_i\rangle = 0\}\,.\end{align*}
Then, for any subset $I\subset I_0$
\begin{equation}\label{char.2}
\beta_I:=\delta_{\gp}-\sum_{i\in I_{-}} \alpha_i-\sum_{i\in I}\alpha_i\,,
\end{equation}
is a highest weight of the spin representation of $K$ with minimal norm.

\ni So there are exactly $1+\# I_0$ such highest weights.
\end{prop}
\proof Let $I$ be a subset of $I_0$. First $\beta_I$ is a weight
of the spin representation of $K$ by \eqref{weights.spin}.
If $\beta_I$ is not a highest weight, then there exists a $K$-positive root $\theta$ such that $\beta_I +\theta$ is a weight.
So there exists a subset $J\subset \{1,\ldots,m\}$ such that
$$\beta_I+\theta =\delta_{\gp}-\sum_{i\in J}\alpha_i\,.
$$
Then,
$$-\sum_{i\in I_{-}\backslash J}\alpha_i-\sum_{i\in I\backslash J}\alpha_i +\theta= -\sum_{i\in J\backslash(I_{-}\cup I)}
\alpha_i\,.$$
Now, it is well known that $\langle \theta, \delta_K\rangle >0$, (cf. for instance \S~10.2 in \cite{Hum}), hence
$\langle -\sum_{i\in I_{-}\backslash J}\alpha_i-\sum_{i\in I\backslash J}\alpha_i +\theta, \delta_K\rangle >0$, whereas
$\langle -\sum_{i\in J\backslash(I_{-}\cup I)} \alpha_i,\delta_K\rangle \leq 0$, hence a contradiction.

\ni By the lemma~\ref{Parth}, there exists a $w_I\in W$ such that
$$\beta_I=w_I\cdot\delta_G-\delta_K= \delta_{\gp}-\sum_{i\in I_{-}} \alpha_i-\sum_{i\in I}\alpha_i\,.$$
Now, using the $W_{G}$-invariance of the scalar product,
\begin{align}\label{sq.norm}\|\beta_I\|^2&=\|\delta_G\|^2+\|\delta_K\|^2-2\langle w_I\cdot\delta_G, \delta_K\rangle\nonumber\\
&=\|\delta_G\|^2+\|\delta_K\|^2-2\,\langle\beta_I+\delta_K, \delta_K\rangle\nonumber\\
&=\|\delta_G\|^2+\|\delta_K\|^2 -2\,\langle \delta_G,\delta_K\rangle +2\,\sum_{i\in I_{-}}\langle
\alpha_i,\delta_K\rangle\nonumber\\
&= \|\delta_{\gp}\|^2+2\,\sum_{i\in I_{-}}\langle\alpha_i,\delta_K\rangle\,.
\end{align}
Hence all the highest weights $\beta_I$, with $I\subset I_0$ have the same norm.
In order to prove that among the highest weights of the spin representation of $K$, they are those with
lower norm, we use the same argument as in \cite{Mil06}.

\ni Let $\theta_1,\ldots,\theta_p$ be an enumeration of the $K$-positive roots. Let
$w\in W$ (or $\in W_G$ as well). By the lemma~3.1 in \cite{Mil06}, using the expression of $w$ in reduced form,
one has
$$w\cdot \delta_G=\delta_G-\sum_{i\in I_w} \alpha_i-\sum_{j\in J_w} \theta_{j}\,,$$
where $I_w$ is a subset of $\{1,\ldots\,m\}$, and $J_w$ a subset of $\{1,\ldots\,p\}$.

\ni Since $\langle \theta_j,\delta_K\rangle>0$, $j=1,\ldots,p$, as we remark before, and since
$\langle \alpha_i,\delta_K\rangle >0$ if $i$ does not belong to $I_{-}$ or $I_0$, one gets
\begin{align}\label{ineqs}
\langle w\cdot\delta_G, \delta_K\rangle &\leq \langle \delta_G-\sum_{i\in
I_w\cap(I_{-}\cup I_0)}\alpha_i,\delta_K\rangle\\
&\leq \langle\delta_G-\sum_{i\in
I_{-}\cup I_0}\alpha_i,\delta_K\rangle\nonumber\\
&\leq \langle\delta_G-\sum_{i\in
I_{-}}\alpha_i,\delta_K\rangle=\langle w_{\emptyset}\cdot\delta_G,\delta_K\rangle\nonumber\,.
\end{align}
Hence
\begin{align*}\| w\cdot\delta_G- \delta_K\|^2&=\|\delta_G\|^2+\|\delta_K\|^2-2\,\langle w\cdot\delta_G,
\delta_K\rangle\\
&\geq \|\delta_G\|^2+\|\delta_K\|^2-2\,\langle w_{\emptyset}\cdot\delta_G,
\delta_K\rangle=\|\beta_{\emptyset}\|^2\,.
\end{align*}
Note that if the above inequality is an equality, then all the inequalities \eqref{ineqs} are equalities, hence
$J_w=\emptyset$ and $I_w=I_{-}\cup I$, where $I\subset I_0$, so $w\cdot\delta_G-\delta_K=\beta_I$.

\ni As the above result is valid for any $w\in W_G$, it may be concluded that for any subset $I\subset I_0$,
 $$\min_{w\in W_G}\| w\cdot\delta_G- \delta_K\|^2=\min_{w\in W}\| w\cdot\delta_G- \delta_K\|^2=\|\beta_I\|^2\,.$$

\qed

\ni Now the proof of \eqref{low.eig} consists in the following steps (\cite{Mil05}):
\begin{enumerate}
 \item
 For any $w\in W$ such that the highest weight $w\cdot \delta_G-\delta_K$ has minimal norm, $\delta_G-w^{-1}\cdot\delta_{K}$ is the
dominant weight of an irreducible representation $\gamma$  of $G$.
 \item The restriction of $\gamma$ to $K$ contains in its decomposition into irreducible components a representation with
dominant weight $w\cdot \delta_G-\delta_K$.
 \item  The Casimir eigenvalue for $\gamma$ is given by $c_{\gamma}=2\, \| w\cdot\delta_G- \delta_K\|^2+\frac{n}{16}$
(hence equal to $2\, \|\delta_{\gp}\|^2+4\, \sum_{i\in I_{-}}\langle\alpha_i,\delta_K\rangle +\frac{n}{16}$, by \eqref{sq.norm}).
 \item  The above Casimir eigenvalue gives the lowest eigenvalue of the Dirac operator.
\end{enumerate}
In the preparation of this paper, we found a gap in the proof of the first item.
We give a different proof in appendix, which is indeed based on the result of proposition~\ref{high.w.char}.

\subsection{The space of eigenvectors of the Dirac operator corresponding to the lowest eigenvalue}

\ni In order to understand the action of a form characterized by the holonomy on the eigenvectors of the Dirac operator
for the lowest eigenvalue, we now review some well-known results (see \cite {BH3M} for details).

\ni First, recall that a spinor field $\Psi$ on $G/K$ may be viewed as a function
$$\Psi : G\longmapsto \Sigma\,,\qquad \forall g\in G\,,\; \forall k\in K\,,\quad
 \Psi(gk)=\widetilde{\rho}(k^{-1})\cdot \Psi(g)\,.$$

 \ni Denoting by $\Sigma_w$ the irreducible $K$-space of $\Sigma$ with dominant weight $\beta_w$, $w\in W$, and by
$\Pi_w$ the projection $\Sigma\rightarrow\Sigma_w$, any spinor field $\Psi$ decomposes into
\be\label{decomp}
\Psi =\sum_{w\in W} \Psi_w\,,\qquad \Psi_w:=\Pi_w\circ \Psi\,.\ee

\ni Since the restricted holonomy group of $G/K$ is $K$, the spin Levi-Civita connection $\nabla$ respects the above decomposition,
hence by the Lichnerowicz-Schr\"odinger formula,
$$\D^2=\nabla^{*}\nabla+\frac{\mathrm{Scal}}{4}\,,$$
 if a spinor field $\Psi$ is an eigenspinor of the Dirac operator $\D$ for the eigenvalue $\lambda$, then
each non-trivial component $\Psi_w$ in the above decomposition is an eigenspinor of $\D^2$ for the eigenvalue $\lambda^2$.

\ni Now, as it was recalled above (see \eqref{eigenv}), any eigenvalue $\lambda$ of the Dirac operator $\D$ corresponds (up to equivalence) to
an irreducible complex $G$-representation $\gamma: G\rightarrow \mathrm{GL}_{\mathbb{C}}(V_\gamma)$, whose restriction $\mathrm{Res}_K^G(\gamma)$ to the subgroup $K$,
contains in its decomposition into irreducible parts, a representation with dominant weight $\beta_w$, $w\in W$.
The corresponding eigenvectors are given by a pair $(v_\gamma, A_\gamma)$, where $v_\gamma\in V_\gamma$ and
$A_\gamma\in \mathrm{Hom}_K(V_\gamma, \Sigma):=\{A\in \mathrm{Hom}_{\bC}(V_\gamma, \Sigma)\,,\; \forall k\in K\;,\; A\circ\gamma(k)=\widetilde{\rho_K}(k)\circ A\}$, giving
raise to the spinor field
$$\Psi_{v_\gamma,A_\gamma}: G\longrightarrow \Sigma\,,\qquad \Psi_{v_\gamma,A_\gamma}(g)=A_\gamma\big(\gamma(g^{-1})\cdot v_\gamma\big)\,.$$
Moreover
$$\dim\,\mathrm{Hom}_K(V_\gamma, \Sigma)=\sum_{w\in W} \mathrm{mult.}(\mathrm{Res}_K^G(\gamma),\widetilde{\rho_{K}}_{|w})\,,$$
where $\mathrm{mult.}(\mathrm{Res}_K^G(\gamma),\widetilde{\rho_{K}}_{|w})$ denotes the multiplicity of the irreducible
representation $\widetilde{\rho_{K}}_{|w}$ in $\mathrm{Res}_K^G(\gamma)$.

\ni So, if $\mathrm{mult.}(\mathrm{Res}_K^G(\gamma),\widetilde{\rho_{K}}_{|w})\not=0$, then
the component $\Psi_w$ in the decomposition of $\Psi$ given by \eqref{decomp} is a non-trivial eigenvector of $\D^2$ for
$\lambda^2$.

\ni All that discussion applies to the irreducible $G$-representation $\gamma$ with dominant weight $\delta_G-w^{-1}\cdot\delta_K$, where $\beta_w=w\cdot\delta_G-\delta_K$
 is a highest weight of the spin representation of $K$ with minimal norm (cf. prop.~\ref{prooflem}). As we recall it above, this irreducible $G$-representation gives raise to the lowest
eigenvalue $\lambda$ of the Dirac operator. Moreover by prop.~\ref{proofrem}, the representation $\mathrm{Res}_K^G(\gamma)$ contains in its decomposition into irreducible components
all the irreducible $K$-representations with dominant weights \eqref{char.2}.

\ni Hence we may conclude
\begin{lem}\label{lem.Psi_I} For any subset $ I\subset I_0 =\{i\in\{1,\ldots,n\}\,,\, \langle \delta_K,\alpha_i\rangle = 0\}$,
denoting by $\Sigma_I$ the irreducible component of $\Sigma$ with highest weight $\beta_I=\delta_{\gp}-\sum_{i\in I_{-}} \alpha_i-\sum_{i\in I}\alpha_i$,
there exists a spinor field $\Psi_I:G\rightarrow \Sigma_I$, such that $$\D^2 \Psi_I=\lambda^2\, \Psi_I\,,$$
where $\lambda$ is the lowest eigenvalue of the Dirac operator.
\end{lem}

\ni Now let $\Omega$ be a parallel form on $G/K$. First, $\Omega$ may be viewed as a $K$-equivariant function
$$\Omega : G\longrightarrow \Lambda^{*}(\gp)\,,\qquad \forall g\in G\,,\forall k\in K\,,\quad \Omega(gk)=\alpha(k^{-1})\cdot\Omega(g),.$$

\ni By the fundamental principle of holonomy, parallel forms correspond to $K$-invariants of $\Lambda^{*}(\gp)$. Hence there exists a $K$-invariant
$\mathbf{\Omega}$ in $\Lambda^{*}(\gp)$ such that $\Omega$ is the constant function
$$\forall g\in G\,,\qquad \Omega(g)=\mathbf{\Omega}\,.$$
Note that as $\mathbf{\Omega}$ is $K$-invariant, $\Omega$ is $K$-equivariant since
$$\forall g\in G\,,\;\forall k\in K\,,\quad \Omega(gk)=\mathbf{\Omega}=\alpha(k^{-1})\cdot\mathbf{\Omega}=\alpha(k^{-1})\cdot\Omega(g)\,.$$

\ni The form $\Omega$ acts on a spinor field $\Psi$ by Clifford multiplication, giving a spinor field $\Omega\cdot\Psi$ defined by the function
$$\Omega\cdot\Psi :G\longrightarrow\Sigma\,,\qquad \Omega\cdot\Psi (g):=\Omega(g)\cdot\Psi(g)=\mathbf{\Omega}\cdot\Psi(g)\,,$$
where $\mathbf{\Omega}$ is viewed as an element of the Clifford algebra, which acts on the spinor $\Psi(g)$ by means of the standard representation of this algebra.

\ni Since $\mathbf{\Omega}$ is $K$-invariant, the Schur lemma implies that the action of $\Omega$ on spinor fields with values in
some $K$-irreducible subspace $\Sigma_w$, $w\in W$, of $\Sigma$, is a scalar multiple of the identity, the value of the scalar depending only of $w\in W$.

\ni In particular, the spinor fields $\Psi_I$, $I\subset I_0$, introduced in lemma~\ref{lem.Psi_I}, which are eigenspinors for the square of the Dirac operator for the lowest eigenvalue,
are also eigenvectors for the action of any parallel form.

\section{Proof of the result}
\subsection{The K\"ahler case}

On a spin K\"ahler manifold of complex dimension $m$, the K\"ahler $2$-form $\Omega$ acts fiberwise on spinors as an anti-hermitian operator with eigenvalues
$i\,(m-2r)$, $r=0,\ldots,m$, \cite{Kir86}. Hence
\be \label{ineq.K}\min_{\Psi\not=0}\frac{\|\Omega\cdot \Psi\|^2}{\|\Psi\|^2}=\begin{cases}0\,,\quad \text{if $m$ is even}\\
1\,,\quad \text{if $m$ is odd}\,.\end{cases}\ee

\ni We are going to prove that \eqref{ineq.K} is verified for one of the spinor fields $\Psi_I$ of lemma~\ref{lem.Psi_I}.

\ni First, an irreducible symmetric space $G/K$ is K\"ahler if and only if $K$ has a center $Z\simeq\mathrm{U}_1$, \cite{KN2}.

\ni Let $\gp_{\bC}^{+}$, (resp. $\gp_{\bC}^{-}$) be the space generated by root-vectors corresponding to the positive non compact roots
(resp. negative non compact roots). Any element $H$ in the Lie algebra $\mathfrak{z}$ of the center has a $K$-invariant adjoint action on $\gp_{\bC}^{+}$,
(resp. $\gp_{\bC}^{-}$), hence by the Schur lemma acts as a scalar multiple of identity.
The element $H$ is chosen such that $\mathrm{ad}(H)_{|\gp_{\bC}^{+}}= i\, \mathrm{id}$
(hence $\alpha_j(H)=i$, $j=1 ,\ldots, m$) and $\mathrm{ad}(H)_{|\gp_{\bC}^{-}}= -i\, \mathrm{id}$.
This action defines a $K$-invariant homomorphism $\mathbf{J}$ of $\gp$ such that $\mathbf{J}^2=-\mathrm{id}$, which induces a K\"ahler structure on $G/K$.
The K\"ahler form is then defined by the $K$-invariant $\mathbf{\Omega}$ corresponding to $\alpha_{*}(H)$ by the isomorphism
$\Lambda^2(\gp)\simeq \mathfrak{so}(\gp)$.
Hence the action of the K\"ahler form on spinor fields is given by the action of $\mathbf{\Omega}$ on $\Sigma$,
which corresponds to $2$ times the action of $H$ on $\Sigma$ by the spinor representation of $K$, since viewed as a $2$-form,
$\mathbf{\Omega}$ is identified with an element of the Clifford algebra, whereas viewed as $\alpha_{*}(H)\in \mathfrak{so}(\gp)$, it acts
on spinors by the isomorphism $\xi_{*}:\mathfrak{spin}(\gp)\rightarrow \mathfrak{so}(\gp)$, which generates a factor $2$.
Finally, as $H$ belongs to the Lie algebra of the maximal torus\footnote{since the center of $K$ is the intersection of the maximal tori.} $T$,
 the K\"ahler form acts on spinor fields with values in $\Sigma_w$, $w\in W$, as a scalar multiple of identity,
the eigenvalue being given by $2\, \beta_w(H)$.

\ni Hence we only have to prove that there exists a subset $I\subset I_0=\{i\in\{1,\ldots,m\}\,,$ $\, \langle \delta_K,\alpha_i\rangle = 0\}$, such that
$\beta_I(H)=0$ if $m$ is even and $\beta_I(H)=\pm i/2$, if $m$ is odd, where $\beta_I=\delta_{\gp}-\sum_{i\in I_{-}} \alpha_i-\sum_{i\in I}\alpha_i$.

\ni Let $I_{+}=\{i\in\{1,\ldots,m\}\;,\; \langle \alpha_i,\delta_K\rangle>0\}$. Then
$$\beta_I(H)=\frac{1}{2}\, \sum_{j\in I_{+}} \alpha_j(H)-\frac{1}{2}\, \sum_{j\in I_{-}} \alpha_j(H)
+\frac{1}{2}\, \sum_{j\in I_{0}\backslash I} \alpha_j(H)-\frac{1}{2}\, \sum_{j\in I} \alpha_j(H)\,.$$
\begin{lem} The sets $I_{+}$ and $I_{-}$ have the same number of elements.
\end{lem}
\proof Let $\Phi_{K}^{-}$ be the set of negative roots of $K$. There exists an element $w_0$ in the Weyl group of $K$ (hence in the Weyl group of $G$)
which sends $\Phi_{K}^{+}$ to $\Phi_{K}^{-}$, see for instance theorem~3.1.9. in \cite{GW09}.

\ni Note that, as $H$ belongs to the Lie algebra $\mathfrak{z}$ of the center $Z$ of $K$, one has for any $K$-root $\theta$,
$\theta(H)=0$, since for any root vector $Y_{\theta}$,
$$0=[H,Y_{\theta}]=\theta(H)\, Y_{\theta}\,.$$
From this remark, one deduces that $w_0 (\Phi_{\gp}^{+})\subset \Phi_{\gp}^{+}$. Indeed, as $H$ belongs to $\mathfrak{z}$,
if $k\in K$ is some representative of $w_0$, one has since $\mathrm{Ad}(k^{-1})\cdot H=H$, and $\alpha_j(H)=i$, $j=1,\ldots,m$,
$$w_0\cdot \alpha_j(H)=\alpha_j(\mathrm{Ad}(k^{-1})\cdot H)=\alpha_j(H)=i\,,\quad j=1,\ldots,m\,.$$

\ni Now let $j\in I_{+}$ so that $\langle \alpha_j,\delta_K\rangle >0$. Then
$\langle w_0\cdot\alpha_j,w_0\cdot \delta_K\rangle>0$, so as  $w_0\cdot \delta_K=-\delta_K$, one gets
$$\langle w_0\cdot\alpha_j, \delta_K\rangle <0\,,$$
hence there exists $i_j\in I_{-}$ such that $\alpha_{i_j}=w_0\cdot\alpha_j$. This defines a one-to-one correspondence between
$I_{+}$ and $I_{-}$.

\qed

\ni Hence since $\alpha_j(H)=i$, $j=1,\ldots,m$, we obtain from the lemma
$$\beta_I(H)=\frac{1}{2}\, \sum_{j\in I_{0}\backslash I} \alpha_j(H)-\frac{1}{2}\, \sum_{j\in I} \alpha_j(H)\,.$$

\ni Now if $m$ is even, then by the result of the lemma, the set $I_0$ has an even number of elements.
If $I_0=\emptyset$, then $\beta_{\emptyset}(H)=0$. If $I_0\not =\emptyset$, then choosing a subset $I\subset
I_0$ such that $\#\, I=\frac{1}{2}\,\#\, I_0$, one gets $\beta_I(H)=0$.

\ni If $m$ is odd, then by the result of the lemma, the set $I_0$ has an odd number $2r+1$ of elements.
Choosing now a subset $I\subset I_0$ such that $\#\, I=r$, (resp. $r+1$) , one gets
$\beta_I(H)=\frac{1}{2}\, i$, (resp. $-\frac{1}{2}\, i$),  and the result is proved.

\subsection{The Quaternion-K\"ahler case}

\ni A Quaternion-K\"ahler manifold  is a  $n=4m$-dimensional Riemannian manifold $(M,g)$ whose restricted holonomy group is contained in the group
$\mathrm{Sp}_m.\mathrm{Sp}_1=\mathrm{Sp}_m\times_{\mathbb{Z}_2}\mathrm{Sp}_1$, $m\geq 2$. This group is identified with a
subgroup of $\mathrm{SO}_{4m}$ by the representation
$$(A,q)\in \mathrm{Sp}_m.\mathrm{Sp}_1 \longmapsto \Big(x\in \bH^m\simeq \bR^{4m}  \mapsto A x \bar{q}\Big)\,.$$

\ni Let $\mathbf{i},\mathbf{j},\mathbf{k}$ be the standard basis of imaginary quaternions. The action on the right of
$-\mathbf{i},-\mathbf{j},-\mathbf{k}$ on $\bH^m$ defines three hermitian operators $\mathbf{I}$, $
\mathbf{J}$, $\mathbf{K}$, verifying the same multiplication rules as the imaginary quaternions.
 The space $\boldsymbol{\mathcal{Q}}$ generated by $\mathbf{I}$, $
\mathbf{J}$, $\mathbf{K}$ is $K$-invariant, hence by transport on the fibres, it defines a globally parallel subbundle
$\mathcal{Q}(M)$ of the bundle $\mathrm{End}(TM)$. By transporting the operators $\mathbf{I}$, $
\mathbf{J}$, $\mathbf{K}$ on fibres with the help of a trivialization, one gets three local
almost complex structures $I$, $J$, $K$, for which the metric $g$ is hermitian,
verifying the same multiplication rules as the imaginary quaternions. Using the metric, one obtains three local
$2$-forms $\Omega_I$, $\Omega_J$, $\Omega_K$. Now, the $4$-form
$$\Omega=\Omega_I\wedge\Omega_I+\Omega_J\wedge\Omega_J+\Omega_K\wedge\Omega_K\,,$$
is well-defined over $M$, parallel and non-degenerate, \cite{Kr66}, \cite{Bon}.

\ni On any spin Quaternion-K\"ahler manifold this $4$-form $\Omega$
acts fiberwise on spinors as an hermitian operator with eigenvalues
$6m-4r(r+2)$, $r=0,\ldots,m$, \cite{HM95a}. Hence
\be \label{ineq.QK}\min_{\Psi\not=0}\frac{\|(\Omega-6m\,\mathrm{id})\cdot \Psi\|^2}{\|\Psi\|^2}=0\,.\ee

\ni We are going to prove that \eqref{ineq.QK} is verified for one of the spinor fields $\Psi_I$ of lemma~\ref{lem.Psi_I}.

\ni Compact symmetric spaces with a Quaternion-K\"ahler structure were classified by J. A. Wolf in \cite{Wol65}.
It is well known in the theory of representations of compact groups that any root associated to the choice of a maximal torus
gives raise to a subgroup of $G$ isomorphic to $\mathrm{Sp}_1$. J. A. Wolf has shown that
 compact symmetric spaces with a Quaternion-K\"ahler are all inner symmetric space of type I of the form $G/K$, where $G$ is
a simple group and $K=K_1\,\mathrm{Sp}_1$, where $K_1$ is the centralizer of $\mathrm{Sp}_1$ in $G$.
The subgroup $\mathrm{Sp}_1$ of $K$ in
consideration here being defined by the maximal root $\beta$ (for a fixed ordering of roots).

\ni Indeed, let $H_\beta\in \gt$ such that for any $H\in \gt$, $\langle H_\beta, H\rangle =-i\,\beta(H)$.
Then $\|H_\beta\|^2=-i\, \beta(H_\beta)=\|\beta\|^2$. Let $H_{\beta}^{\circ}:=2/{\|\beta\|^2}\, H_{\beta}$.

\ni Let $X_{\beta}$ be a root-vector for the root $\beta$. There exists a root-vector $X_{-\beta}$ for the root
$-\beta$ such that $[X_\beta,X_{-\beta}]=-i\, H_{\beta}^{\circ}$.

\ni Then $(H_{\beta}^{\circ}, Y_{\beta}:=i\, (X_{-\beta}+X_{\beta}),
Z_{\beta}:= X_{-\beta}-X_{\beta})$ defined a basis of a subagebra of $\gg$ isomorphic to $\mathfrak{sp}_1$ as
$$[H_{\beta}^{\circ},Y_{\beta}]=2\, Z_{\beta}\,,\quad [H_{\beta}^{\circ},Z_{\beta}]=-2\, Y_{\beta}\quad\text{and}
\quad [Y_{\beta},Z_{\beta}]=2\,H_{\beta}^{\circ}\,.$$
Now, the condition that $\beta$ is the maximal root implies
that  $\mathrm{ad} (H_\beta^{\circ})_{|\gp_{\bC}^{+}}=i\, \mathrm{id}$, \cite{Wol65}, so
the action of $H_{\beta}^{\circ},Y_{\beta}$ and $Z_{\beta}$ on $\gp$ by $\alpha_{*}:\gk\rightarrow
\mathfrak{so}(\gp)$, induces three hermitian operators $\mathbf{I}$, $
\mathbf{J}$, $\mathbf{K}$, verifying the same multiplication rules as the vectors $\mathbf{i},\mathbf{j},\mathbf{k}$
of the standard basis of imaginary quaternions.
The space $\boldsymbol{\mathcal{Q}}$ generated by $\mathbf{I}$, $
\mathbf{J}$, $\mathbf{K}$, which is $K$-invariant, generates the Quaternion-K\"ahler structure on $G/K$.

\ni Identifying $\mathbf{I}$, $\mathbf{J}$ and $\mathbf{K}$ with $2$-forms
 $\mathbf{\Omega_I}$, $\mathbf{\Omega_J}$, $\mathbf{\Omega_K}$,
via the metric, one gets the $K$-invariant $4$-form on $\gp$
$$\mathbf{\Omega}=\mathbf{\Omega_I}\wedge\mathbf{\Omega_I}+\mathbf{\Omega_J}\wedge\mathbf{\Omega_J}+
\mathbf{\Omega_K}\wedge\mathbf{\Omega_K}\,,$$
which induces the Quaternion-K\"ahler parallel fundamental $4$-form $\Omega$ on $G/K$.

\ni Now, if the symmetric space has a spin structure, the $4$-form $\mathbf{\Omega}$ acts on the spinor space
$\Sigma$ as the operator (\cite{HM95b})
$$\mathbf{\Omega}=6m\, \mathrm{id}+\mathbf{\Omega_I}\cdot\mathbf{\Omega_I}+
\mathbf{\Omega_J}\cdot\mathbf{\Omega_J}+\mathbf{\Omega_K}\cdot\mathbf{\Omega_K}\,,$$
where the $2$-forms $\mathbf{\Omega_I}$, $\mathbf{\Omega_J}$, $\mathbf{\Omega_K}$ act by Clifford multiplication.
Hence \footnote{here again the presence of
the scalar factor $4$ is due to the use of the isomorphism $\xi_{*}:\mathfrak{spin}(\gp)\rightarrow
\mathfrak{so}(\gp)$, when the two-forms $\mathbf{\Omega_I}$, $\mathbf{\Omega_J}$, $\mathbf{\Omega_K}$ are
identified with $H_{\beta}^{\circ}, Y_{\beta},Z_{\beta}$,
acting on spinors by the representation $\widetilde{\rho}_{*}$ of $\gk$.}
\be\label{Omega.QK}
\mathbf{\Omega}-6m\, \mathrm{id} =4\,\Big(\widetilde{\rho}_{*}(H_{\beta}^{0})^2
+\widetilde{\rho}_{*}(Y_{\beta})^2+
\widetilde{\rho}_{*}(Z_{\beta})^2\Big)
 \,.\ee

\ni Note that the second term in the r.h.s. of the above equation is the Casimir operator\footnote{up to some normalization.}
of the representation $\widetilde{\rho}_{*}$ restricted to $\mathfrak{sp}_1$.

\ni Expressing the r.h.s. of \eqref{Omega.QK} in the basis $(H_{\beta}^{\mathbf{c}}:=-i\,H_{\beta}^{\circ}, X_{\beta},X_{-\beta})$
 of $\mathfrak{sp}_{1,\bC}\simeq\mathfrak{sl}_{2,\bC}$, one gets
\begin{align}\label{Cas2}\mathbf{\Omega}-6m\, \mathrm{id}&=-4\,\Big(\widetilde{\rho}_{*}(H_{\beta}^{\mathbf{c}})^2
+2\, \widetilde{\rho}_{*}(X_{\beta})\circ\widetilde{\rho}_{*}(X_{-\beta})+2\,
\widetilde{\rho}_{*}(X_{-\beta})\circ\widetilde{\rho}_{*}(X_{\beta})\Big)\nonumber\\
&=-4\,\Big(\widetilde{\rho}_{*}(H_{\beta}^{\mathbf{c}})^2
+2\,\widetilde{\rho}_{*}(H_{\beta}^{\mathbf{c}}) +4\,
\widetilde{\rho}_{*}(X_{-\beta})\circ\widetilde{\rho}_{*}(X_{\beta})\Big)\,.
\end{align}
Hence we may conclude
\begin{lem}\label{OmegaPsiI} For any subset $I\subset I_{0}$,
$$(\Omega-6m\, \mathrm{id})\cdot \Psi_{I}=0\Longleftrightarrow \beta_{I}(H_\beta^{\circ})=0\,.$$
\end{lem}
\proof By the Schur lemma, $\mathbf{\Omega}-6m$ acts on the $K$-irreducible space $\Sigma_I$ as a scalar
multiple of identity. If $c_I$ is the eigenvalue, one then has $(\Omega-6m\, \mathrm{id})\cdot \Psi_{I}=
c_I \, \Psi_I$. To compute the eigenvalue, one applies \eqref{Cas2} to a highest weight vector of $\Sigma_{I}$.
Since the action of $\widetilde{\rho}_{*}(X_{\beta})$ is zero on such a vector,
 whereas $\widetilde{\rho}_{*}(H_{\beta}^{\mathbf{c}})$ acts by a non-negative integer
multiple of identity on it, one has $c_I=0$ if and only if
$\widetilde{\rho}_{*}(H_{\beta}^{\mathbf{c}})$ acts trivially, hence the result.

\qed

\ni Let $I_{+}=\{i\in\{1,\ldots,m\}\;,\; \langle \alpha_i,\delta_K\rangle>0\}$. Then
\be\label{betaI}
\beta_I(H_\beta^{\circ})=\frac{1}{2}\, \sum_{j\in I_{+}} \alpha_j(H_\beta^{\circ})-
\frac{1}{2}\, \sum_{j\in I_{-}} \alpha_j(H_\beta^{\circ})
+\frac{1}{2}\, \sum_{j\in I_{0}\backslash I} \alpha_j(H_\beta^{\circ})
-\frac{1}{2}\, \sum_{j\in I} \alpha_j(H_\beta^{\circ})\,.
\ee

\begin{lem}\label{lemQK} Apart from $G_2/SO_4$, for any Quaternion-K\"ahler compact spin symmetric space, one has
 $$\#\, I_{-}+\#\, I_{0}=\#\,I_{+}\,.$$
\end{lem}
\proof Note first that
\be\label{orth.beta} \forall \theta \in \Phi_{K}\,,\; \theta\not =\pm\beta\,,
\quad \langle \beta ,\theta \rangle =0\,. \ee
Indeed if $X_\theta$ is a root-vector for the root $\theta$, one has $[H_\beta, X_\theta]=0$, since $K_1$ is the centralizer
of $\mathrm{Sp}_1$ in $K$. So $\theta(H_\beta)=0$. Now, let $H_\theta\in \gt$ be such that for any $H\in \gt$
$\langle H_\theta, H\rangle= -i\, \theta(H)$. We then have $\langle H_\theta,H_\beta\rangle=0$, hence
$\langle \theta ,\beta \rangle =0$.

\ni By this remark, positive non compact roots $\alpha$ are then characterized by the condition
 $\langle \alpha, \beta\rangle =\frac{1}{2}\, \|\beta\|^2$.

\ni Let $\alpha_i$ be a positive non compact root. Then $\beta-\alpha_i$ is a positive compact root. It is a root
since $\langle \alpha_i, \beta\rangle >0$ (cf. for instance\footnote{or note that
$\sigma_\beta(\alpha_i)=\alpha_i-\beta$, where $\sigma_\beta$ the reflection across the hyperplane $\beta^{\bot}$.}
\S~9.4 in \cite{Hum}). And it is positive since $\langle \beta-\alpha_i,\beta\rangle=\frac{1}{2}\, \|\beta\|^2$.
Note furthermore that $\beta-\alpha_i\not=\alpha_i$, since otherwise $2\alpha_i=\beta$ should be a root, which is
impossible.

\ni Now by \eqref{orth.beta}, $\langle \delta_K, \beta\rangle= \frac{1}{2}\, \|\beta\|^2$, hence
$$\langle \delta_K, \beta-\alpha_i\rangle=\frac{1}{2}\, \|\beta\|^2-\langle \delta_K,\alpha_i\rangle\,.$$

\ni Hence, if $j\in I_{-}\cup I_{0}$, then $\langle \delta_K, \beta-\alpha_j\rangle\geq  \frac{1}{2}\, \|\beta\|^2$,
hence $\beta-\alpha_j=\alpha_{i_j}$, with $i_j\in I_{+}$. We thus get an injective map $I_{-}\cup I_{0}\rightarrow
I_{+}$, so we may conclude $\#\, I_{-}+\#\, I_{0}\leq\#\,I_{+}$.
\ni On the other hand, if $\langle \delta_K,\alpha_j\rangle >0$, then as $\delta_K$ is an integral
weight\footnote{since $\delta_K$ is a difference of integral weights:
$\delta_K=w\cdot\delta_G-(w\cdot\delta_G-\delta_K)$, $w\in W$.}, one has
$\langle \delta_K,\alpha_j\rangle\geq\frac{1}{2}\, \|\alpha_j\|^2$, hence
$$\langle \delta_K, \beta-\alpha_j\rangle\leq \frac{1}{2}\, (\|\beta\|^2-\|\alpha_j\|^2)\,.$$
Now as $G$ is a simple group, the root system is irreducible and there are at most two root lengths (see for instance
\S~10.4. in \cite{Hum}).

\ni If all the roots have same length, then $\langle \delta_K, \beta-\alpha_j\rangle \leq 0$, so
$\beta-\alpha_j=\alpha_{i_j}$, where $i_j\in I_{-}\cup I_{0}$. In this case, there is an injective
map $I_{+}\rightarrow I_{-}\cup I_{0}$, so $\#\,I_{+}\leq\#\, I_{-}+\#\, I_{0}$, and the result is proved.

\ni Now by the result \cite{Wol65}, using furthermore the result \cite{CG},
 the list of spin compact Quaternion-K\"ahler symmetric spaces is given by
 \vspace{5mm}
\begin{center}
\begin{tabular} {|c|c|c|c|c|} \hline
 $G$&$K$&$ G/K$& $\dim\, G/K$&Spin structure  \\  \hline $\mathrm{Sp}_{m+1}$
&$\mathrm{Sp}_{m}\times\mathrm{Sp}_{1}$& Quaternionic & $4m\,
(m\geq 1)$&Yes (unique) \\ &&projective&&\\ &&space
$\mathbb{H}P^{m}$&&\\ \hline
$\mathrm{SU}_{m+2}$&$S(\mathrm{U}_{m}\times\mathrm{U}_{2})$ &
Grassmannian &$4m\, (m\geq 1)$&iff $m$ even
\\
&&$\mathrm{Gr}_{2}(\mathbb{C}^{m+2})$&&unique in that case\\
\hline$\mathrm{Spin}_{m+4}$& $\mathrm{Spin}_{m}\mathrm{Spin}_{4}$&
Grassmannian &$4m\, (m\geq 3)$ &iff $m$ even,
\\
&&$\widetilde{\mathrm{Gr}}_{4}(\mathbb{R}^{m+4})$&&unique in that
case\\ \hline $\mathrm{G}_{2}$&$\mathrm{SO}_{4}$&&$8$&Yes
(unique)\\&&&&\\ \hline $\mathrm{F}_{4}$ &
$\mathrm{Sp}_{3}\mathrm{SU}_{2}$&&$28$&No\\&&&&\\ \hline $
\mathrm{E}_{6}$& $\mathrm{SU}_{6}\mathrm{SU}_{2}$&&$40$&Yes
(unique)\\&&&&\\ \hline $\mathrm{E}_{7}$&
$\mathrm{Spin}_{12}\mathrm{SU}_{2}$&&$64$&Yes (unique)\\ &&&&\\
\hline $\mathrm{E}_{8}$ & $\mathrm{E}_{7}\mathrm{SU}_{2}$&&$112$&
Yes (unique)\\&&&&\\ \hline
\end{tabular}\end{center}
\vspace{5mm}
\ni Note that all of them are inner as it was noticed in \cite{Wol65}.

\ni Now, apart from\footnote{there are two root lengths for $G=\mathrm{F}_4$, but the
corresponding symmetric space is not spin.}
$\mathrm{Sp}(m+1)$ and $\mathrm{G}_2$,
 there are only one root length for the groups $G$
in the above list, hence the result is proved for the corresponding symmetric spaces.

\ni So it remains to prove the result for quaternionic projective spaces $\mathbb{H}P^m=$\\
$\mathrm{Sp}(m+1)/\mathrm{Sp}(m)\times \mathrm{Sp}_1$, $m\geq 1$.

\ni We consider the standard maximal torus $T$ of $\mathrm{Sp}_{m+1}$ made up of diagonal matrices
with entries of the form $\mathrm{e}^{\beta\, \mathbf{i}}:=\cos(\beta)+\sin(\beta)\,\mathbf{i}$, $\beta\in\bR$.
We denote by $(x_0,x_1,\ldots, x_m)$ the standard basis of $\gt^{*}$ such that the value of $x_k$ on
a diagonal matrix with entries $(\beta_{0}\, \mathbf{i},\ldots, \beta_{m}\, \mathbf{i})$ is $\beta_k$,
$k=0,\ldots, m$. We set $\widehat{x}_k=i\, x_k$. The scalar product on $i\, \gt^{*}$
induced by the Killing form sign-changed verifies $\langle
\hat{x}_i,\hat{x}_j\rangle=\frac{1}{4(m+2)}\,\delta_{ij}$.
We choose as positive roots
$$\widehat{x}_i\pm\widehat{x}_j\,,\quad 0\leq i<j\leq m\,; \qquad 2\, \widehat{x}_i\,,\quad
0\leq i\leq m\,.$$
The roots $\widehat{x}_i-\widehat{x}_{i+1}$, $0\leq i\leq m-1$, and $2\,\hat{x}_m$ then define a basis of the root
system.

\ni In order to avoid a re-ordering of roots, we consider $K=\mathrm{Sp}_1\times \mathrm{Sp}_m$, (instead
of $\mathrm{Sp}_m\times \mathrm{Sp}_1$) in such a way that the positive compact roots are
$$\widehat{x}_i\pm\widehat{x}_j\,,\quad 1\leq i<j\leq m\,; \qquad 2\, \widehat{x}_i\,,\quad
0\leq i\leq m\,.$$

\ni Then, the positive non compact roots are
$$\widehat{x}_0\pm\widehat{x}_k\,,\quad 1\leq k\leq m\,.$$
The maximal root is $\beta=2\, \hat{x}_0$. Note that $\widehat{x}_0\pm\widehat{x}_k=\beta -
(\widehat{x}_0\mp\widehat{x}_k)$, and $\langle \beta, \widehat{x}_0\pm\widehat{x}_k\rangle=1/2\, \|\beta\|^2$.

\ni Now $\delta_K=\hat{x}_0+\sum_{k=1}^{m} (m-k+1)\, \hat{x}_k$. Hence it is easy to verify that
\begin{align*}
&\langle \delta,\widehat{x}_0+\widehat{x}_k\rangle >0\,,\qquad 1\leq m\,,\\
&\langle \delta,\widehat{x}_0+\widehat{x}_k\rangle <0\,,\qquad 1\leq m-1\,,\\
&\langle \delta,\widehat{x}_0+\widehat{x}_m\rangle =0\,.
\end{align*}
So $\#\, I_{-}+\#\, I_{0}=\#\,I_{+}$.

\qed

\ni Going back to \eqref{betaI}, one deduces from the above lemma, as
$\mathrm{ad} (H_\beta^{\circ})_{|\gp_{\bC}^{+}}=i\, \mathrm{id}$,
$$\beta_{I_{0}}(H_{\beta}^{\circ})=\frac{1}{2}\,i\, (\#\,I_{+}-\#\, I_{-}-\#\, I_{0})=0\,.$$
Hence by the lemma~\ref{OmegaPsiI}, $(\Omega-6m\, id)\cdot \Psi_{I_{0}}=0$, and the result
of proposition~\ref{result} is proved.

\subsection {The case of $\mathrm{G}_2/\mathrm{SO}_4$}
The group $\mathrm{SO}_4$ is identified with $\mathrm{Sp}_1\cdot\mathrm{Sp}_1\simeq
\mathrm{Sp}_1\times_{\mathbb{Z}_2}\mathrm{Sp}_1$.
The inclusion $\mathrm{Sp}_1\cdot\mathrm{Sp}_1\subset \mathrm{Sp}_2\cdot\mathrm{Sp}_1$ is not
the ``{natural}'' one since the group acts irreducibly on $\bH^2=\bR^8$ with highest weight
$3\, \omega_1+\omega_2$, where $(\omega_1,\omega_2)$ is the standard basis of fundamental weights
corresponding to the half spinors representations, see for instance \S 11 in \cite{Sal89}.

\ni We use the result of \cite{See99}, where all the roots data for $G$ and $K$ are expressed
in the basis $(\omega_1,\omega_2)$ by
\begin{gather*}\Phi_G^{+}=\{2\,\omega_1, -3\, \omega_1+\omega_2,-\omega_1+\omega_2,\omega_1+\omega_2,2\,\omega_2,
3\, \omega_1+\omega_2\}\\
\Phi_K^{+}=\{2\,\omega_1,2\,\omega_2\}\,,\\
\Phi_{\gp}^{+}=\{-3\, \omega_1+\omega_2, -\omega_1+\omega_2,\omega_1+\omega_2, 3\,\omega_1+\omega_2\}\,.
\end{gather*}

\ni The scalar product on weights induced by the Killing form sign-changed is given by
$$\langle a_1\,\omega_1+a_2\,\omega_2, b_1\,\omega_1+b_2\,\omega_2\rangle=
\frac{1}{16}\,\left(\frac{1}{3}\,a_1b_1+a_2b_2\right)\,.$$

\ni The highest weights of the spin representation are obtained by means of the Partha-saraty formula,
\cite{See99}
$$4\,\omega_1\,,\quad 2\, \omega_2\,\quad \text{and}\quad 3\,\omega_1+\omega_2\,.$$
The half spin representation $\Sigma_{8}^{-}$ is irreducible with highest weight $3\,\omega_1+\omega_2$,
the half spin representation $\Sigma_{8}^{+}$ decomposes into two components with respective highest
weights $4\, \omega_1$ and $2\,\omega_2$. So, denoting by $H$ the standard representation of
$\mathrm{Sp}_1$, the spin representation decomposes as
\be\label{spinG2}\Sigma_{8}=(\mathrm{S}^3\, H\otimes \mathrm{S}^1\, H)\oplus (\mathrm{S}^4\, H\otimes
\mathrm{S}^{0}\, H)\oplus (\mathrm{S}^0\, H\otimes \mathrm{S}^2\,H)\,.\ee

\ni There are two weights for which the norm is minimal: $\beta_{\emptyset}=2\, \omega_2=\delta_{\gp}$ and
$\beta_{I_{0}}=3\,\omega_1+\omega_2=\delta_{\gp}-(-3\,\omega_1+\omega_2)$.

\ni Hence the spinor field $\Psi_{\emptyset}$ verifies $\mathcal{D}\Psi_{\emptyset}= \lambda^2\, \Psi_{\emptyset}$, where $\lambda$
is the lowest eigenvalue of the Dirac operator, and is also a section of the bundle corresponding to the component
$\mathrm{S}^0\, H\otimes \mathrm{S}^2\,H$ in the decomposition \eqref{spinG2}, on which the action (on the first component)
of $\mathrm{Sp}_1$ is trivial, thus our ``{holonomy criterion}'' is verified. However, it has to be noticed that
the Quaternion-K\"ahler structure of $\mathrm{G}_2/\mathrm{SO}_4$ corresponds to the action of $\mathrm{Sp}_1$ on the
second component, and that this criterion is not verified for that action.

\ni Indeed, it easy to see that $\theta_1=2\,\omega_1$ and $\theta_2=-3\omega_1+\omega_2$ are simple $G$-roots, and that
the maximal root is $\beta:=2\, \omega_2=3\, \theta_1+2\,\theta_2$. Now note that

$$\begin{cases} \langle \beta, \alpha\rangle =0 &\quad \text{if $\alpha\in \Phi_K$, $\alpha\not =\pm \beta$,}\\
\langle \beta, \alpha\rangle =\frac{1}{2}\, \|\beta\|^2&\quad \text{if $\alpha\in \Phi_{\gp}^{+}$,}
\end{cases}$$
hence we are in the description of Wolf spaces given above. It may be checked that for any positive non compact
positive root $\alpha$, $\langle \delta_K,\alpha\rangle=0$, $1/24$, $1/12$ or $1/8$, so
$\#\,I_{+}=3$, $\#\,I_{0}=1$ and $\#\,I_{-}=0$, hence the result of lemma~\ref{lemQK} is not verified here.

\ni Moreover note that, by definition of $H_{\beta}^{\circ}$,
 $\beta_{\emptyset}(H_{\beta}^{\circ})=\frac{2\, i}{\|\beta\|^2}\,
\langle \beta, \beta_{\emptyset}\rangle\not =0$ and $\beta_{I_{0}}(H_{\beta}^{\circ})=\frac{2\, i}{\|\beta\|^2}\,
\langle \beta, \beta_{I_0}\rangle\not =0$, so by lemma~\ref{OmegaPsiI}, the criterion is not verified
for the action of $\mathrm{Sp}_1$ on the second component.

\ni As an additional argument, one may remark that since $\langle \beta, 4\,\omega_1 \rangle=0$,
$4\, \omega_1(H_{\beta}^{\circ})=0$, so by the proof of lemma~\ref{OmegaPsiI},
the action of $\mathbf{\Omega} -6m\, \mathrm{id}$ on spinors is zero only on the $K$ irreducible subspace of $\Sigma$
with highest weight $4\,\omega_1$. But, by the results of \cite{See99}, this is not a highest weight for the
restriction to $K$ of the $\mathrm{G}_2$ irreducible
representation giving raise to the lowest eigenvalue of the Dirac operator.

\ni Indeed, the action of $\mathrm{Sp}_1$ on the first component may also be described in terms of
the action of fundamental $4$-form on spinors by Clifford multiplication. This is a consequence of a more general
result.

\begin{lem}\label{CasimirK}
The eigenvalue of the Casimir operator of the spin representation of $K$  has the same value for all
 irreducible components.
\end{lem}
\proof By the Parthasaraty formula, any highest weight of the spin reresentation of $K$ has the form
$\beta_w:=w\cdot\delta_G-\delta_K$, where $w$ is an element of the Weyl group $W_G$ of $G$. By the
Freudenthal formula, the eigenvalue of the Casimir operator of $K$ acting on the $K$-representation
with highest weight $\beta_w$ is given by
$$\langle \beta_w+2\, \delta_K, \beta_w\rangle\,.$$
Hence, by the $W_G$-invariance of the scalar product,
\begin{align*}
\langle \beta_w+2\, \delta_K, \beta_w\rangle&=\langle w\cdot\delta_G+\delta_K,
w\cdot\delta_G-\delta_K\rangle\\
&= \|\delta_G\|^2-\|\delta_K\|^2\,.
\end{align*}

\qed

\ni Hence, since $\|\delta_G\|^2-\|\delta_K\|^2=\frac{1}{2}$ here,
the Casimir operator $\mathcal{C}_K$ of the spin representation of $K$ acts on each irreducible component
as $\frac{1}{2}\,\mathrm{id}$.

\ni Denote by $\mathcal{C}_1$ (resp. $\mathcal{C}_2$) the Casimir operator of the restriction
of the spinor representation to the first
(resp. second) $\mathrm{Sp}_1$ component in $K$, (for the scalar product given by
the Killing form sign-changed of $\mathrm{G}_2$). One has
$\mathcal{C}_K=\mathcal{C}_1+\mathcal{C}_2$. The eigenvalue of $\mathcal{C}_1$
($\mathcal{C}_2$) acting on $\mathrm{S}^k H$ is a scalar multiple of $k(k+2)$ times the
identity. Using the result of the above lemma, one gets by \eqref{spinG2}
$${\mathcal{C}_1}_{|\mathrm{S}^{k}H}=\frac{1}{48}\, k(k+2)\, \mathrm{id}\quad \text{and}
\quad {\mathcal{C}_2}_{|\mathrm{S}^{k}H}=\frac{1}{16}\, k(k+2)\,\mathrm{id}\,.$$
Hence using \eqref{Cas2}, one obtains
$$\mathcal{C}_1= \frac{1}{2}\,\mathrm{id}-\mathcal{C}_2
= \frac{1}{64}\, (\mathbf{\Omega}+20\, \mathrm{id})\,.
$$
So for the symmetric space $\mathrm{G}_2/\mathrm{SO}_4$, the holonomy criterion is valid
for the $4$-fundamental form acting (by Clifford multiplication) on spinors as
$\Omega+20\, \mathrm{id}$.

\subsection{The Cayley plane $\mathrm{F}_4/\mathrm{Spin}_9$.}

\ni The above ``{holonomy criterion}'' is also valid for the Cayley plane $\mathrm{F}_4/\mathrm{Spin}_9$.
 As it is said in \cite{Bes}, the $\mathrm{Spin}(9)$ holonomy is extremely special
since a Riemannian manifold whose holonomy group is contained in $\mathrm{Spin}(9)$ is either
 flat or (locally) isometric to the Cayley plane
$\mathrm{F}_4/\mathrm{Spin}_9$, or its non-compact dual,
 \cite{Ale}, \cite{BG72}. There is an analogy between
 $\mathrm{Spin}(9)$-structures (for $16$-dimensional manifolds) and Quaternion-K\"ahler structures,
since such a structure on a manifold $M$ may be characterized by the existence of a $9$-dimensional subbundle of the bundle
$\mathrm{End}(TM)$, with local sections $I_{\alpha}$, $1\leq \alpha\leq 9$, satisfying
$$I_{\alpha}^2=\mathrm{id}\quad,\quad I_{\alpha}^{*}=I_{\alpha}\quad,\quad
I_{\alpha}I_{\beta}=-I_{\beta}I_{\alpha}\;,\; \alpha\not =\beta\,,$$
\cite{Fri01}.
There is also an analogy by the existence of a canonical $8$-form $\Omega$, which corresponds
to the unique parallel $8$-form on $\mathrm{F}_4/\mathrm{Spin}_9$. This canonical $8$-form was first
introduced in \cite{BG72} by means an integral formula. Explicit algebraic expressions of this form
 are far from being simple \cite{BPT85}, \cite{AM96}, \cite{CLMGM}, \cite{PP12}. Roughly speaking, the form is
 constructed by means of the K\"ahler $2$-forms associated to the almost complex structures
$J_{\alpha \beta}:=I_{\alpha}\circ I_{\beta}$, \cite{CLMGM}, \cite{PP12}.
To avoid an explicit expression, we will use here the fact that it may be expressed in terms of ``{higher
Casimir operators}'', see \S 125 and \S 126 in \cite{Zel} or \cite{Hom04}.

\ni Indeed, this parallel $8$-form is induced by an $\mathrm{Ad}(K)$-invariant $8$-form $\mathbf{\Omega}$ on
$\gp=\bR^{16}$. From the expression of $\Omega$ given in \cite{CLMGM}, \cite{PP12},
the action of $\mathbf{\Omega}$ on spinors is (up to a shift by a scalar multiple of the identity) a sum of terms of the form
$ \widetilde{\rho_K}_{*}(\Omega_{\alpha_1\beta_1})\circ \widetilde{\rho_K}_{*}(\Omega_{\alpha_2\beta_2})\circ
\widetilde{\rho_K}_{*}(\Omega_{\alpha_3\beta_3})
\circ\widetilde{\rho_K}_{*}(\Omega_{\alpha_4\beta_4})$, where the $\Omega_{\alpha\beta}$ are the K\"ahler $2$-forms, identified
with elements of the Lie algebra $\mathfrak{spin}_9$, associated to the almost complex structures $J_{\alpha
\beta}$. We thus may identify the action of $\mathbf{\Omega}$ on spinors as the action of an element
of the universal enveloping algebra $\mathrm{U}(\mathfrak{spin}_{9,\bC})$, also denoted $\mathbf{\Omega}$.
Furthermore this element belongs to the center
$\mathfrak{Z}$ of $\mathrm{U}(\mathfrak{spin}_{9,\bC})$ since $\mathbf{\Omega}$ is
$\mathrm{Ad}(K)$-invariant.

\ni Now it is known that the center $\mathfrak{Z}$ is algebraically generated by a system of
 ``{higher Casimir elements}'', cf. \S.125 and \S.126 in \cite{Zel}, which we briefly introduce in
that context.

\ni Let $(e_i)_{1\leq i\leq 9}$ be the standard basis of $\bR^9$. Denote by $e_{ij}=e_i\wedge e_j$ the standard
basis of $\mathfrak{spin}_9\simeq\mathfrak{so}_9\simeq \Lambda^2(\bR^9)$. For any non negative integer $q$, consider
$e_{ij}^q\in \mathrm{U}(\mathfrak{so}_{9,\bC})$ defined by
$$e_{ij}^q:=\begin{cases}\sum_{1\leq k_1,\ldots,k_{q-1}\leq 9} e_{i k_1}e_{k_1 k_2}\ldots e_{k_{q-1} j}\,,& q\geq 1\,,\\
\delta_{ij}\,, & q=0\,.
\end{cases}
$$

\ni Then the trace of $e_{ij}^q$, $C_q:=\sum_{i=1}^{9} e_{ii}^q$ belongs to the center $\mathfrak{Z}$. For $q=0$,
$C_0=9$, for $q=1$, $C_1=0$ and for $q=2$, $C_2=\sum_{ij} e_{ij}e_{ji}$ is the usual Casimir element.

\ni It may be shown that the center $\mathfrak{Z}$ is algebraically generated by $C_2$, $C_4$, $C_6$ and $C_8$, \cite {Zel}.

\ni Now, the above description of the element $\mathbf{\Omega}\in \mathrm{U}(\mathfrak{spin}_{9,\bC})$ shows that
it is expressed only in terms of $C_2$ and $C_4$. If we consider an irreducible component of the spin representation, then
the Schur lemma implies that $\widetilde{\rho_K}_{*}(C_2)$ and $\widetilde{\rho_K}_{*}(C_4)$ acts on it as scalar multiples of identity.
Furthermore, we already note
in lemma~\ref{CasimirK}, that $\rho_{*}(C_2)$  acts on each irreducible component as a scalar multiple of identity with the same eigenvalue.
Finally, to prove our holonomy criterion, we only have to examine the eigenvalues of $\rho_{*}(C_4)$ on the irreducible components
of the spin representation of $K$. Those eigenvalues may be computed with the help of a formula given in \cite{CGH00}.

\ni Denoting by $(e_1,e_2,e_3, e_4)$ the standard basis of $\bR^4$,
the root system of $\mathrm{F}_4$ is the set of elements $x=\sum_{i=1}^{4} x_i\, e_i$ with integer or half-integer
coordinates in $\bR^4$ such that $\|x\|^2=1$ or $2$, \cite{Hum},\cite{BMP85}. Using the results in \cite{CG},
we may consider
\begin{gather*}\Phi_G^{+}=\{e_i\,,\, i=1,2,3,4\,;\,
e_i\pm e_j\,,\, 1\leq i<j\leq 4\,;\,\frac{1}{2}\, (e_1\pm e_2\pm e_3\pm e_4)\}\\
\Phi_K^{+}=\{e_i\,,\, i=1,2,3,4\,;\,
e_i\pm e_j\,,\, 1\leq i<j\leq 4\}\,,\\
\Phi_{\gp}^{+}=\{\frac{1}{2}\, (e_1\pm e_2\pm e_3\pm e_4)\}\,.
\end{gather*}
Thus
\begin{gather*}\delta_G=\frac{1}{2}\, (11\,e_1+5\, e_2+3\, e_3+e_4)\,,\\
\delta_K=\frac{1}{2}\, (7\,e_1+5\, e_2+3\, e_3+e_4)\,,\\
\delta_{\gp}=2\, e_1\,.
\end{gather*}
The scalar product induced by the Killing form sign-changed is a scalar multiple
of the restriction to the set of roots of the usual scalar product on $\bR^4$.
Using the ``{strange}'' formula of Freudenthal and De Vries, \cite{FdV}, one obtains
$\|\delta_G\|^2=\frac{\dim \gg}{24}=13/6$, hence this scalar product is given by
$$\langle \sum_{i=1}^{4} x_i\,e_i,\sum_{i=1}^{4} y_i\,e_i\rangle= \frac{1}{18}\,\sum_{i=1}^{4} x_i\, y_i\,.$$
With the help of the Parthasaraty formula, it is easy to find the highest weights of the spin representation of $K$.
The half spin representation $\Sigma_{16}^{-}$
is irreducible with highest weight $\frac{1}{2}(3\, e_1+e_2+e_3+e_4)$, whereas the half spin representation
 $\Sigma_{16}^{+}$ decomposes into two components with highest weights $e_1+e_2+e_3$ and $2\,e_1$.

\ni Note that the sets $I_{-}$ and $I_{0}$ have both only one element since for any positive non compact root
$\alpha$, one has
$\langle \delta_K, \alpha\rangle <0\Longleftrightarrow \alpha= \frac{1}{2}\,(e_1-e_2-e_3-e_4)$, and
$\langle \delta_K, \alpha\rangle =0\Longleftrightarrow \alpha= \frac{1}{2}\,(e_1-e_2-e_3+e_4)$.

\ni Indeed, there are two weights for which the norm is minimal
\begin{align*}\beta_{\emptyset}&=\frac{1}{2}(3\,
e_1+e_2+e_3+e_4)=\delta_{\gp}-\frac{1}{2}\,(e_1-e_2-e_3-e_4)\,,\\
\intertext{and}
\beta_{I_0}&=e_1+e_2+e_3=\delta_{\gp}-\frac{1}{2}\,(e_1-e_2-e_3-e_4)-\frac{1}{2}\,(e_1-e_2-e_3+e_4)\,.
\end{align*}
By the result of \cite{Mil05}, the square of the first eigenvalue of the Dirac operator is then given by
$2\, \|\beta_{\emptyset}\|^2+2=2\, \|\beta_{I_{0}}\|^2+2$, hence
\begin{prop} On the symmetric space $\mathrm{F}_4/\mathrm{Spin}_9$ endowed with the Riemannian metric induced
by the Killing form of $\mathrm{F}_4$ sign-changed, the square of the first eigenvalue $\lambda$ of the Dirac
operator verifies
$$\lambda^2=\frac{7}{3}=\frac{7}{6}\, \frac{\mathrm{Scal}}{4}\,.$$
\end{prop}

\ni We now verify the holonomy criterion by determining explicitly the eigenvalues of the operator
$\widetilde{\rho_K}_{*}(C_4)$. This may be done with the help of a formula given in \cite{CGH00}
which applies here as follows.

\ni For each dominant weight $\mu$ of $\mathrm{SO}_9$ (relatively to the usual maximal torus), denote
by $d(\mu)$ the dimension of a complex $K$-irreducible representation $(\rho_\mu, V(\mu))$ with highest weight $\mu$.

\ni Let $V(\mu_{\circ})$ be the standard representation of $\mathrm{SO}_9$ corresponding to the
weight $\mu_{\circ}=(1,0,0,0)$. For each highest weight
$\beta_1$, $\beta_2$, $\beta_3$ of the spin representation of $K$ where
$$\beta_1=\frac{1}{2}(3\, e_1+e_2+e_3+e_4)\,,\quad
\beta_2=e_1+e_2+e_3\, \quad \text{and} \quad \beta_3=2\, e_1\,,$$
denote by $\Pi_i$ the set of highest weights occurring in the decomposition into irreducible components
of the tensor product $V(\beta_i)\otimes V(\mu_{\circ})$. So
$$V(\beta_i)\otimes V(\mu_{\circ})=\sum_{\lambda_{ij}\in \Pi_i} V(\lambda_{ij})\,.$$

\ni For each highest weight $\lambda_{ij}$, let $m(\lambda_{ij})$ be the so called conformal
weight given by
$$m(\lambda_{ij}):=\frac{1}{2}\,(9-\|\delta_K+\lambda_{ij}\|^2+\|\delta_K+\beta_i\|^2-1)\,,$$
where $\|\cdot\|$ is the standard norm on the sets of weights of $\mathrm{SO}_9$: if
$\mu=(\mu_1,\mu_2,\mu_3,\mu_4)$, then $\|\mu\|^2=\sum_{i=1}^{4} \mu_i^2$.

\ni The eigenvalue $c_{\beta_i,k}$ of the operator ${(\rho_{\beta_i})}_{*}(C_k)$ is then given by
(see for instance \cite{Hom04})
$$c_{\beta_i,k}=\frac{1}{d(\beta_i)}\, \sum_{\lambda_{ij}\in \Pi_i} \big(-m(\lambda_{ij})\big)^{k}\,
d(\lambda_{ij})\,.$$

\ni The values of $c_{\beta_i,k}$, $i=1,2,3$, $k=2,3,4$ are given in the following table.
They were obtained with the help of the LiE Program\footnote{\url{http://www-math.univ-poitiers.fr/~maavl/LiE/}.}.

\begin{center}
\begin{tabular}{|c|c|c|c|}
\hline
&$\beta_1=\frac{1}{2}(3,1,1,1)$ &$\beta_2=(1,1,1,0)$& $\beta_3=(2,0,0,0)$\\
\hline
$c_{\beta_i,2}$&36&36&36\\
\hline
$c_{\beta_i,3}$&-126&-126&-126\\
\hline
$c_{\beta_i,4}$&1863/2&684&1404\\
\hline
\end{tabular}
\end{center}
The first line is not a surprise by lemma~\ref{CasimirK}. The second line too since $C_3$ may be expressed in
terms of $C_2$. Finally, the last line shows that the lowest eigenvalue of the operator
$\widetilde{\rho_K}_{*}(C_4)$ is obtained for the space with highest weight $\beta_2=(1,1,1,0)$.
But $\beta_2=\beta_{I_0}$, which as we saw it before, is linked to the first eigenvalue of the Dirac operator.
Hence the ``{holonomy criterion}'' of prop.~\ref{result} is verified.

\section{Appendix}
\ni Recall that the highest weights of the spin representation of $K$ are given by
$$\beta_w=w\cdot\delta_G-\delta_K\,,\quad w\in W\,,$$
where $W$ is the subset of the Weyl group $W_G$ defined by
$$W:=\{w\in W_{G}\;;\;
w\cdot \Phi_{G}^{+}\supset \Phi_{K}^{+}\}\,.$$

\ni In this appendix, we review the proof of the following result.

\begin{prop}\label{prooflem}Let $w_{0}\in W$ be such that
\begin{equation}\label{w0}
 \|\beta_{w_{0}}\|^{2}=\min_{w\in W}\|\beta_{w}\|^{2}\,.
\end{equation}
Then the weight
$$\beta_{w_{0}}^{G}:=w_{0}^{-1}\cdot\beta_{w_{0}}=\delta_{G}-w_{0}^{-1}\cdot\delta_{K}\,,$$
is $G$-dominant.
\end{prop}

\ni We will also justify the following remark~:
\begin{prop}\label{proofrem} Let $w_{0}$ and $w_{1}\in W$ be such that
$$\|\beta_{w_{1}}\|^{2}= \|\beta_{w_{0}}\|^{2}=\min_{w\in W}\|\beta_{w}\|^{2}\,.$$
Then the $G$-dominant weights
$$\beta_{w_{0}}^{G}:=w_{0}^{-1}\cdot\beta_{w_{0}}=\delta_{G}-w_{0}^{-1}\cdot\delta_{K}
\quad \text{and}\quad
\beta_{w_{1}}^{G}:=w_{1}^{-1}\cdot\beta_{w_{0}}=\delta_{G}-w_{1}^{-1}\cdot\delta_{K}
$$ verify
$$\beta_{w_{0}}^{G}=\beta_{w_{1}}^{G}\,.$$

\end{prop}

\subsection{Review of the proof of prop.~\ref{prooflem}}

\ni The main error in the proof given in \cite{Mil05} concerns a technical assumption which is used at the end,
asserting that the highest weights with same minimal length may be ordered with  respect to the usual
order of $K$-weights. As it can be seen in \eqref{char.2}, this is actually not correct since
two such highest weights may only differ by a \textit{non compact} positive root\footnote{whereas, if it is a root,
a sum of compact roots is necessarily a compact one.}. Before that, there is also an imprecise statement, asserting that a
certain weight of the spin representation does lie in the $W_G$-orbit of a highest weight, which makes the proof not
satisfactory.

\ni Hence we propose here an alternative proof we think correct.

\proof First note that considering the sets
$$\Lambda_{w_0}^{+}:=\{\theta\in \Phi_{G}^{+}\,,\, w_0\cdot\theta\in \Phi_{\gp}^{+}\}\quad\text{and}\quad
\Lambda_{w_0}^{-}:=\{\theta\in \Phi_{G}^{+}\,,\, -w_0\cdot\theta\in \Phi_{\gp}^{+}\}\,,$$ one has
$$\beta_{w_0}=\frac{1}{2}\,\sum_{\theta\in\Phi_{G}^{+}\backslash w_{0}^{-1}(\Phi_{K}^{+})}w_0\cdot\theta=
\frac{1}{2}\,\sum_{\theta\in\Lambda_{w_0}^{+}}w_0\cdot\theta+
\frac{1}{2}\,\sum_{\theta\in\Lambda_{w_0}^{-}}w_0\cdot\theta\,. $$
Hence since
$$\delta_{\gp}=\frac{1}{2}\,\sum_{\theta\in\Lambda_{w_0}^{+}}w_0\cdot\theta-\frac{1}{2}\,
\sum_{\theta\in\Lambda_{w_0}^{-}}w_0\cdot\theta\,, $$
\begin{equation}\label{char.1}
 \beta_{w_0}=\delta_{\gp}+\sum_{\theta\in\Lambda_{w_0}^{-}}w_0\cdot\theta\,.
\end{equation}

\ni Now by the result of prop.~\ref{high.w.char}, there exists a subset $I_{w_0}\subset I_0$ such that
\be\label{char.1'}
\beta_{w_0}=\delta_{\gp}-\sum_{i\in I_{-}} \alpha_i-\sum_{i\in I_{w_{0}}} \alpha_i\,.
\ee
Comparing \eqref{char.1} and \eqref{char.1'}, we get the following alternative~: if a positive $G$-root is such that
$w_0\cdot\theta$ is a non compact root, then either $w_0\cdot\theta$ is negative and then $w_0\cdot\theta=-\alpha_i$,
$i\in I_{-}\cup I_{w_{0}}$, so $\langle w_0\cdot\theta,\delta_K\rangle \geq 0$, or $w_0\cdot\theta$ is positive and then
 $w_0\cdot\theta=\alpha_i$,
$i\notin I_{-}\cup I_{w_{0}}$, so $\langle w_0\cdot\theta,\delta_K\rangle \geq 0$.

\ni So we may conclude that
\begin{rem}\label{rem} If $\alpha$ is a positive $G$-root such that $w_0\cdot\alpha$ is a non compact root then, whatever
$w_0\cdot\alpha$ is positive or not, one always has
$$\langle \delta_K, w_0\cdot\alpha\rangle\geq 0\,.$$
\end{rem}

\ni  Let $\Pi_{G}=\{\theta_{1},\ldots,\theta_{r}\}\subset \Phi_{G}^{+}$ be the
set of simple roots. It is sufficient to prove that
$2\,\frac{\langle\beta_{w_{0}}^{G},\theta_{i}\rangle}{\langle\theta_{i},\theta_{i}\rangle}$
is a non-negative integer for any simple root $\theta_{i}$. First, as
$T$ is a maximal common torus of $G$ and $K$, $\beta_{w_{0}}$,
which is an integral weight for $K$ is also an integral weight for
$G$. Now since the Weyl group $W_{G}$ permutes the weights,
$\beta_{w_{0}}^{G}=w_{0}^{-1}\cdot\beta_{w_{0}}$ is also a
integral weight for $G$, hence
$2\,\frac{\langle\beta_{w_{0}}^{G},\theta_{i}\rangle}{\langle\theta_{i},\theta_{i}\rangle}$
is an integer for any simple root $\theta_{i}$. So we only have to
prove that $\langle\beta_{w_{0}}^{G},\theta_{i}\rangle=\langle\delta_{G}-w_{0}^{-1}\cdot\delta_{K},\theta_{i}\rangle\geq 0$,
or equivalently (by the $W_{G}$-invariance of the scalar product) that
\begin{equation}\label{ineq1}\langle w_{0}\cdot \delta_{G}-\delta_{K},w_{0}\cdot\theta_{i}\rangle\geq 0\,.
\end{equation}

\ni Let $\theta_{i}$ be a
simple root. Suppose first that $w_{0}\cdot\theta_{i}\in \Phi_{K}$.

\noindent In this case, necessarily
$w_{0}\cdot\theta_{i}\in \Phi_{K}^{+}$, otherwise since $w_{0}\in W$, $-\theta_{i}=w_{0}^{-1}(-w_{0}\cdot\theta_{i})$
should be a positive root.

\noindent Then since $w_{0}\cdot \delta_G-\delta_K$ is $K$-dominant, inequality \eqref{ineq1} is verified in this
case, since $w_{0}\cdot\theta_{i}$ is a linear combination with non-negative integer coefficients of $K$-simple
roots.

\ni Suppose now that $w_{0}\cdot\theta_{i}\notin \Phi_{K}$, that is
$w_{0}\cdot\theta_{i}$ is a non compact root. We are going to prove that
\begin{equation}\label{ineq2}\langle w_{0}\cdot \delta_{G}-\delta_{K},w_{0}\cdot\theta_{i}\rangle<0\,,
\end{equation}
is impossible. Suppose that \eqref{ineq2} is true. Note first that since
$2\,\frac{\langle\delta_{G},\theta_{i}\rangle}{\langle\theta_{i},\theta_{i}\rangle}=1$,
(see for instance \S 10.2 in \cite{Hum}) and since the scalar
product $\langle\cdot ,\cdot\rangle$ is $W_{G}$-invariant, inequation~\eqref{ineq2} is equivalent to
\begin{equation}\label{ineq3}
\langle \delta_{K},w_{0}\cdot\theta_{i}\rangle>\frac{1}{2}\, \|\theta_{i}\|^{2}\,.
\end{equation}
This implies that
$$\langle \delta_{K},w_{0}\cdot\theta_{i}\rangle >0\,.
$$
Now since $\delta_{K}$ is a linear combination with non-negative coefficients of $K$-simple
roots, this implies that there exists a $K$-simple root $\theta'$ such that
\begin{equation}\label{ineq4}\langle\theta',w_{0}\cdot\theta_{i}\rangle >0\,.\end{equation}
By property of roots (cf. for instance
\S 9.4 in \cite{Hum}), this implies that
$$\theta'-w_{0}\cdot\theta_{i}\,,$$
is a root, and moreover a non compact root by the bracket relations $[\gk,\gp]\subset\gp$.

\noindent By the definition of $W$, $\theta:=w_{0}^{-1}\cdot \theta'$ is a positive $G$-root.
 By inequality \eqref{ineq4} $\langle \theta, \theta_i\rangle>0$,
hence $\theta-\theta_i$ is a root, and moreover a positive root.

\noindent Now, the non compact root $\theta'-w_{0}\cdot\theta_{i}$ being the image by $w_0$
of the $G$-positive root $\theta-\theta_i$, the above remark~\ref{rem} applies:
\begin{equation}\label{ineq5}\langle \delta_K,\theta'-w_{0}\cdot\theta_{i}\rangle\geq 0\,.
\end{equation}
Since $\theta'$ is a $K$-simple root, $\langle \delta_K, \theta'\rangle=\frac{1}{2}\,\|\theta'\|^2$, hence we obtain
from inequality \eqref{ineq3}
\begin{equation}\label{ineq6}
 \frac{1}{2}\, \|\theta_i\|^2< \langle \delta_K,w_{0}\cdot\theta_{i}\rangle\leq\frac{1}{2}
\|\theta'\|^2\,.
\end{equation}
Now as $G$ is a simple group, the root system is irreducible and there are at most two root lengths (see for instance
\S~10.4. in \cite{Hum}).

\noindent If all roots are of equal length, then the above inequality is impossible, and the result is proven.

\noindent If there are two distinct root lengths, then the above inequality is only possible if $\theta'$ is a long root
and $\theta_i$ a short one.

\noindent We rewrite inequalities \eqref{ineq6} as
$$1<2\, \frac{\langle \delta_K,w_{0}\cdot\theta_{i}\rangle}{\langle \theta_i,\theta_i\rangle}\leq \frac{\|\theta'\|^2}
{\|\theta_i\|^2}\,.$$
Now $\|\theta'\|^2/\|\theta_i\|^2$ is equal to $2$ or $3$, see for instance \S~9.4 in \cite{Hum}.

\ni Let us examine the case where $\|\theta'\|^2/\|\theta_i\|^2= 2$ first.

\noindent If the inequality \eqref{ineq5} is strict, one obtains
$$1<2\, \frac{\langle \delta_K,w_{0}\cdot\theta_{i}\rangle}{\langle \theta_i,\theta_i\rangle}< \frac{\|\theta'\|^2}
{\|\theta_i\|^2}\,,$$
hence a contradiction, since as $\delta_K=w_0\cdot\delta_G -(w_0\cdot\delta_G-\delta_K)$ is an integral weight,
$2\, \frac{\langle \delta_K,w_{0}\cdot\theta_{i}\rangle}{\langle \theta_i,\theta_i\rangle}$ is an integer.

\noindent Thus \eqref{ineq5} is an equality :
\be\label{ineq5=}
\langle \delta_K,\theta'-w_{0}\cdot\theta_{i}\rangle= 0\,.\ee

\noindent This implies that there exists $j\in I_0$ such that $\theta'-w_{0}\cdot\theta_{i}=\pm \alpha_j$.

\noindent So, by the result of prop.~\ref{high.w.char}
\begin{itemize}
\item either $\theta'-w_{0}\cdot\theta_{i}=\alpha_j$ and $j\in I_{w_0}$ or $\theta'-w_{0}\cdot\theta_{i}=-\alpha_j$
and $j\notin I_{w_0}$ and then
 $$\mu_0:=(w_0\cdot\delta_G-\delta_K)+(\theta'-w_0\cdot\theta_i)\,$$ is a highest weight of
the spin representation with minimal length,
\item either $\theta'-w_{0}\cdot\theta_{i}=\alpha_j$ and $j\notin I_{w_0}$ or $\theta'-w_{0}\cdot\theta_{i}=-\alpha_j$
and $j\in I_{w_0}$ and then
 $$\mu_0:=(w_0\cdot\delta_G-\delta_K)-(\theta'-w_0\cdot\theta_i)\,$$ is a highest weight of
the spin representation with minimal length.
\end{itemize}

\ni In the first case, one gets using $\|\mu_0\|^2=\|w_0\cdot\delta_G-\delta_K\|^2$, \eqref{ineq5=}
and $\langle \delta_G,\theta_i\rangle=\frac{1}{2}\,\|\theta_i\|^2$,
\begin{align*}
0&=2\, \langle w_0\cdot\delta_G-\delta_K,\theta'-w_{0}\cdot\theta_{i}\rangle
+\|\theta'-w_0\cdot\theta_i\|^2\,,\\
&= 2\,\langle w_0\cdot\delta_G,\theta'\rangle-\|\theta_i\|^2+\|\theta'-w_0\cdot\theta_i\|^2\,,\\
&= 2\,\langle \delta_G,w_{0}^ {-1}\cdot\theta'\rangle-\|\theta_i\|^2+\|\theta'-w_0\cdot\theta_i\|^2\,,
 \end{align*}
whereas in the second case, one obtains
$$0= -2\,\langle \delta_G,w_{0}^ {-1}\cdot\theta'\rangle+\|\theta_i\|^2+\|\theta'-w_0\cdot\theta_i\|^2\,,$$

\ni But inequality~\eqref{ineq4} implies as $2\,\frac{\langle \theta',w_0\cdot\theta_i\rangle}{\langle \theta',\theta'\rangle}$
is an integer
$$2\, \langle \theta',w_0\cdot\theta_i\rangle\geq \|\theta'\|^2\,.$$
Hence
$$\|\theta'-w_0\cdot\theta_i\|^2=\|\theta'\|^2-2\, \langle \theta',w_0\cdot\theta_i\rangle+\|\theta_i\|^2\leq \|\theta_i\|^2\,.$$
This implies that $\theta'-w_0\cdot\theta_i$ is a short root, hence
that
$$\|\theta'-w_0\cdot\theta_i\|^2=\|\theta_i\|^2\,,\quad\text{and}\quad 2\, \langle \theta',w_0\cdot\theta_i\rangle=\|\theta'\|^2\,.$$

\ni In the first case, this implies
$$ \langle \delta_G,w_{0}^ {-1}\cdot\theta'\rangle=0\,,$$
which is impossible, since $w_0^{-1}\cdot\theta'$ is a positive root, because $\langle \delta_G,\theta_j\rangle>0$ for any
simple $G$-root $\theta_j$.

\ni In the second case, one obtains
\begin{align*}0&=-2\,\langle w_0\cdot\delta_G,\theta'\rangle
+2\,\|\theta_i\|^2\\
&=-2\,\langle w_0\cdot\delta_G,\theta'\rangle+\|\theta'\|^2\\
&=-2\,\langle w_0\cdot\delta_G,\theta'\rangle+2\, \langle \theta',w_0\cdot\theta_i\rangle\,,
\end{align*}
hence
$$\langle w_0\cdot\delta_G,\theta'\rangle= \langle w_0\cdot\theta_i,\theta'\rangle\,.$$
Let $\sigma_i$ be the reflection across the hyperplane $\theta_i^{\bot}$. Since
$\sigma_i\cdot\delta_G=\delta_G-\theta_i$ (see for instance \S~10.2 in \cite{Hum}), we obtain
$$\langle w_0\sigma_i\cdot\delta_G,\theta'\rangle= 0\,,$$
hence using the $W_G$-invariance of the scalar product
\begin{equation}\label{ineq7}\langle \delta_G, \sigma_i\,w_0^{-1}\cdot\theta'\rangle= 0\,.
 \end{equation}
But $w_0^{-1}\cdot\theta'$ is a positive root which, being a long root, is different from $\theta_i$.
Hence as $\sigma_i$ permutes the positive roots other than $\theta_i$, (see for instance \S~10.2 in \cite{Hum}),
$\sigma_i\,w_0^{-1}\cdot\theta'$ is a positive root. Then as $\langle \delta_G,\theta_j\rangle>0$ for any
simple $G$-root $\theta_j$, \eqref{ineq7} is impossible.

\ni So the result is proven in the case where $\|\theta'\|^2/\|\theta_i\|^2= 2$.

\ni We finally examine the case where $\|\theta'\|^2/\|\theta_i\|^2= 3$. The only simple group for which
this is possible is the group $G_2$, and there is only one symmetric space of type I to be considered:
$G_2/SO_4$. So we may give a direct proof of proposition~\ref{prooflem} in that case.

\ni We use the result of \cite{See99}, where all the roots data for $G$ and $K$ are expressed
in terms of the fundamental weights $(\omega_1,\omega_2)$ of $\mathrm{SO}_4$, corresponding to the half spinors
representations.
\begin{gather}\Phi_G^{+}=\{2\,\omega_1, -3\, \omega_1+\omega_2,-\omega_1+\omega_2,\omega_1+\omega_2,2\,\omega_2,
3\, \omega_1+\omega_2\}\nonumber\\
\Phi_K^{+}=\{2\,\omega_1,2\,\omega_2\}\,.\nonumber\end{gather}
Note that $\theta_1=2\,\omega_1$ and $\theta_2=-3\omega_1+\omega_2$ are simple $G$-roots, and that
$$\delta_G=\omega_1+3\, \omega_2\quad \text{and}\quad \delta_K=\omega_1+\omega_2\,.$$
The scalar product on $i\,\mathfrak{T}^{*}$ induced by the Killing form sign-changed is given by
$$\langle a_1\,\omega_1+a_2\,\omega_2, b_1\,\omega_1+b_2\,\omega_2\rangle=
\frac{1}{16}\,\left(\frac{1}{3}\,a_1b_1+a_2b_2\right)\,.$$
The highest weights of the spin representation are (\cite{See99}):
$$\beta_1=4\,\omega_1\,,\quad \beta_2=2\, \omega_2\,\quad \text{and}\quad \beta_3=3\,\omega_1+\omega_2\,.$$
Note that there are two highest weights for which the norm is minimal: $\beta_2$ and $\beta_3$. One has
$$\beta_2=\delta_G-\delta_K\quad \text{and} \quad \beta_3=\sigma_2\cdot\delta_G-\delta_K\,,$$
where $\sigma_2$ is the reflection across the hyperplane orthogonal to the simple root $\theta_2$.

\noindent Note that
$$\sigma_2^{-1}\cdot\beta_3=\sigma_2\cdot\beta_3=2\, \omega_2\,.$$
Now
$$2\, \frac{\langle 2\,\omega_2, \theta_1\rangle}{\langle\theta_1,\theta_1\rangle}=0\quad \text{and}\quad
2\, \frac{\langle 2\,\omega_2, \theta_2\rangle}{\langle\theta_2,\theta_2\rangle}=1\,,$$
hence the weight $2\,\omega_2=\beta_2=\sigma_2^{-1}\cdot\beta_3$ is $G$-dominant, so the result
is also proven in that case.

\qed

\subsection{Proof of prop.~\ref{proofrem}}
By the result of prop.~\ref{high.w.char}, there exist two distinct elements $w_{0}$ and $w_{1}$ in $W$ such
that $\|\beta_{w_{1}}\|^{2}= \|\beta_{w_{0}}\|^{2}=\min_{w\in W}\|\beta_{w}\|^{2}$, only if $I_{0}\not= \emptyset$.
Hence we suppose $I_{0}\not= \emptyset$. Let $i\in I_0$ and let $I=\{i\}$.

\ni By \eqref{char.2}, $\beta_I=\beta_{\emptyset}-\alpha_I$.

\ni Since $\|\beta_I\|^2=\|\beta_{\emptyset}\|^2$, one obtains
$$2\,\langle \beta_{\emptyset},\alpha_{i}\rangle =\langle \alpha_i,\alpha_i\rangle\,,$$
hence denoting by $\sigma_{\alpha_{i}}$ the reflection across the hyperplane $\alpha_i^{\bot}$,
$$\sigma_{\alpha_{i}}\cdot\beta_{\emptyset}=\beta_{\emptyset}-\alpha_{i}=\beta_{I}\,.$$
On the other hand, as $i\in I_0$, $\langle \delta_K,\alpha_i\rangle=0$, so
$$\sigma_{\alpha_{i}}\cdot \delta_K=\delta_K\,.$$
Hence, considering $w_{\emptyset}$ and $w_I\in W$ such that
$\beta_{\emptyset}=w_{\emptyset}\cdot\delta_G-\delta_K$,
and $\beta_{I}=w_{I}\cdot\delta_G-\delta_K$, one gets
$$w_{I}\cdot\delta_G-\delta_K=\beta_I=\beta_{\emptyset}-\alpha_I=
\sigma_{\alpha_{i}}\cdot\beta_{\emptyset}=\sigma_{\alpha_{i}}w_{\emptyset}\cdot\delta_G-\delta_K\,.$$
This implies
$$w_{I}^{-1}\sigma_{\alpha_{i}}w_{\emptyset}\cdot\delta_G=\delta_G\,.$$
But this is only possible if
$$w_{I}=\sigma_{\alpha_{i}}w_{\emptyset}\,,$$
(see for instance \S~122 in \cite{Zel}). Thus
$$\beta_{I}^{G}=\delta_G-w_{I}^{-1}\cdot\delta_K=\delta_G-w_{\emptyset}^{-1}\sigma_{\alpha_{i}}\cdot\delta_K
=\delta_G-w_{\emptyset}^{-1}\cdot\delta_K=\beta_{\emptyset}^G\,.$$
Repeating the argument, the result follows by induction on the cardinal of $I\subset I_0$.

\qed

\end{document}